\documentclass[11pt]{article}
\usepackage{mathrsfs}
\usepackage{amsfonts}
\usepackage{amssymb,amsmath,amsthm}
\newcommand{\R}{\mathbb{R}}

\newcommand{\N}{\mathbb{N}}
\newcommand{\beq}{\begin{equation}}
\newcommand{\eeq}{\end{equation}}
\newcommand{\bp}{\begin{proof}}
\newcommand{\ep}{\end{proof}}
\newcommand{\bo}{\begin{proposition}}
\newcommand{\eo}{\end{proposition}}
\newcommand{\bl}{\begin{lemma}}
\newcommand{\el}{\end{lemma}}

\newtheorem{theorem}{Theorem}[section]
\newtheorem{corollary}{Corollary}

\newtheorem{lemma}[theorem]{Lemma}
\newtheorem{proposition}{Proposition}

\theoremstyle{definition}

\newtheorem{remark}{Remark}

\begin{document}
\title{\LARGE\bf{Uniqueness and nondegeneracy of ground states
for $(-\Delta)^su+u=2(I_2\star u^2)u$ in $\mathbb{R}^N$ when $s$ is close to 1}$\thanks{{\small This work was partially supported by NSFC(11901532).}}$ }
\date{}
 \author{ Huxiao Luo$\thanks{{\small Corresponding author. E-mail: luohuxiao@zjnu.edu.cn (H. Luo).}}$\\
\small Department of Mathematics, Zhejiang Normal University, Jinhua, Zhejiang, 321004, P. R. China
}
\maketitle
\begin{center}
\begin{minipage}{13cm}
\par
\small  {\bf Abstract:} In this article, we study the uniqueness and nondegeneracy of ground states to the fractional Choquard equation:
$$(-\Delta)^su+u=2(I_2\star u^2)u,\quad x\in\R^N,$$
where $N\in\{3,4,5\}$, $s\in(0,1)$ is sufficiently close to $1$. Our method is to make a continuation argument with respect
to the power $s\in(0,1)$ appearing in $(-\Delta)^s$. This approach is based on [M. M. Fall and E. Valdinoci, Comm. Math. Phys., 329 (2014) 383-404].


\end{minipage}
\end{center}

 {\section{Introduction}}
 \setcounter{equation}{0}
In this article, we study the uniqueness and nondegeneracy of ground states for the fractional Choquard equation
\begin{equation}\label{eq:20220420-e1}
\left\{
\begin{array}{ll}
\aligned
&(-\Delta)^s u+ u=2\left(I_2\star u^{2}\right)u \quad \text{in}~\R^N,\\
&u\in H^s(\R^N),
\endaligned
\end{array}
\right.  \tag{P}
\end{equation}
where $2<N<2+4s$, $s\in(0,1)$ sufficiently close to $1$, $\star$ represents the convolution operation in $\R^N$. $I_\alpha$ is the Riesz potential, which is defined for each $x\in\R^N\setminus\{0\}$ by
$$I_\alpha(\cdot)=  \frac{A_\alpha}{|\cdot|^{N-\alpha}},\quad \alpha\in(0,N),$$
where the normalisation constant $A_\alpha$ ensures that the semigroup property
$$I_\alpha \star I_\beta = I_{\alpha+\beta},\quad  \forall \alpha, \beta> 0~\text{such ~that} ~\alpha+\beta < N,$$
see \cite{MR3625092}.
In addition, $I_2$ is the Green function of the Laplacian $-\Delta$ on $\R^N$. The fractional Laplacian $(-\Delta)^s$ defined as
$$(-\Delta)^su(x) =C_{N,s} P.V. \int_{\R^N}\frac{u(x)- u(y)}{|x - y|^{N+2s}}dy$$
where P.V. denotes the Cauchy principal value,
$$C_{N, s} =\left(\int_{\R^N}\frac{1-cos\zeta_1}{|\zeta|^{N+2s}}d\zeta\right)^{-1}.$$
 The fractional Laplacian $(-\Delta)^s$ can be viewed as a pseudo-differential
operator of symbol $|\xi|^{2s}$,
\beq\label{eq1:20221118}
\widehat{(-\Delta)^su}(\xi)=|\xi|^{2s}\widehat{u}(\xi),\quad \xi\in\R^N,
\eeq
where $\widehat{u}$ is the Fourier transform of the function $u$, see \cite{CS,CGM,GZ,DPV} for details.

When $s=1$, $N=3$, \eqref{eq:20220420-e1} changes to be the Schr\"{o}dinger-Newton equation
  \beq\label{20220609-VU3}
-\Delta u+u=2u\left(I_2\star u^2\right)~~\hbox{in}~\R^3.
\eeq
\eqref{20220609-VU3} describe standing waves for the Hartree equation or the coupling of the Schr\"{o}dinger's equation under a classical Newtonian gravitational potential.
Pekar \cite{Pekar} used \eqref{20220609-VU3} to describe the quantum theory of a polaron at rest. Later, Choquard \cite{Lieb1976/77} used it to model one-component plasmas.

The uniqueness of the solution for \eqref{20220609-VU3} has been solved, see in \cite{Lieb1976/77,Ma2010,MR1677740,MR2492602}. In \cite{WY}, by developing Lieb's methods \cite{Lieb1976/77}, T. Wang and T. S. Yi established uniqueness of the positive radial solutions to
\beq\label{20220913-e0}
-\Delta u+u=\left(\frac{1}{|x|^{N-2}}\star u^2\right)u~~\hbox{in}~\R^N,~~u\in H^1(\R^N).
\eeq
This, together with L. Ma and L. Zhao's classification results \cite{Ma2010}, shows the uniqueness of the positive solutions when $N=3,4,5.$

The nondegeneracy of the solution for (\ref{20220609-VU3}) was proved by Lenzmann \cite[Theorem 1.4]{Lenzmann2009}, see also Tod-Moroz \cite{MR1677740} and Wei-Winter \cite{MR2492602}. This nondegeneracy result has
been generalized to \eqref{20220913-e0}
for general dimensions $N=3,4,5$, see \cite{MR4204770}.

In \cite{Xiang}, C.-L. Xiang obtained uniqueness and nondegeneracy results for ground states of the Choquard
equation
$$-\Delta u+u = \left(\frac{1}{|x|}\star |u|^p\right) |u|^{p-2}u~~\hbox{in}~\R^3,$$
provided that $p > 2$ and $p$ is sufficiently
close to $2$.

As far as we know, there is no uniqueness and nondegeneracy result for fractional Choquard equations at present, except for the pseudorelativistic Hartree equation studied in \cite{Lenzmann2009}.

The uniqueness and nondegeneracy of ground state for the following fractional Schr\"{o}dinger equations have been solved,
\begin{equation}\label{eq:bec}
(-\Delta)^su +u = u^p,
\end{equation} where $s\in (0, 1)$ and $p\in(1,2^*_s-1)$.
In \cite{FV}, M. M. Fall and E. Valdinoci proved that if $s$ is sufficiently close to $1$ equation \eqref{eq:bec} possesses a unique minimizer, which is nondegenerate.
Soon after, R. L. Frank, E. Lenzmann and L. Silvestre \cite{FLS} proved that the uniqueness and nondegeneracy of ground states for \eqref{eq:bec} hold for any $s \in (0, 1)$.

Inspired by \cite{FV}, we study the uniqueness and nondegeneracy of ground states for \eqref{eq:20220420-e1}
 by the continuation argument with respect to the power $s\in(0,1)$.  Compared with \cite{FV},  there are some new difficulties.
 \begin{itemize}
\item[(1)] Firstly, the non-local property of Hartree term prevents us from getting the uniform boundedness of the minimizer. To overcome this problem, we rewrite the equation \eqref{eq:20220420-e1} into the system \eqref{20221223-e0}, and obtain uniform boundedness by using the double blow-up method;
\item[(2)] The core problem in our setting is the proof of the nondegeneracy for ground states. Due to the presence of the Hartree term, we can't use the Perron-Frobenius-type arguments in \cite{FV} directly.
To overcome this difficulty, we rely on more complex analysis for the action of the linearized system \eqref{20221219-e1} with respect to decomposition into spherical harmonics, which is the main innovation of this article, see Lemma \ref{3.9};
\item[(3)] After decomposing the solution of the linearized system \eqref{20221219-e1} into spherical harmonics, we use the nondegeneracy of the solution of the limit equation ($s\nearrow 1$) to prove that the radial symmetric part $w_s \equiv 0$. To this end, we need to show that $\|w_s\|_{H^s(\R^N)}$ is bounded. Due to the nonlocality of Hartree term, we can't prove the boundedness similar to \cite{FV}. However, we overcome this difficulty by the semigroup property of the Riesz potential, see Sec. 3.3.
\end{itemize}

We point out that the existence of ground states for \eqref{eq:20220420-e1} has been obtained. In \cite{PGM}, P. d'Avenia, G. Siciliano and M. Squassina studied the ground state for the general fractional Choquard equation
\beq\label{eq1:2022914}
(-\Delta)^su + \omega u = (I_\alpha\star |u|^p)|u|^{p-2}u,\quad u \in H^s(\R^N),
\eeq
where $\omega> 0$, $N \geq 3$, $\alpha\in(0,N)$, $p > 1$ and $s \in (0, 1)$. By assuming that
\beq\label{eq2:2022914}
1+\frac{\alpha}{N}<p<\frac{N+\alpha}{N-2s},
\eeq
the authors obtained the following results.
\begin{proposition}\label{th1.1}  (\cite[Theorem 1.1]{PGM})
\begin{itemize}
\item[(i)]
There exists a ground state $u \in H^s(\R^N)$ to problem $(\ref{eq1:2022914})$  which is
positive, radially symmetric and decreasing;
\item[(ii)] $u \in L^1(\R^N)$ and moreover if $s \leq 1/2$, $u \in C^{0,\sigma}(\R^N)$ for some $\sigma\in
(0, 2s)$, if $s > 1/2$, $u \in C^{1,\sigma}(\R^N)$ for some $\sigma\in (0, 2s- 1)$;
\item[(iii)] If $p \geq 2$, there exists $C > 0$ such that
$u(x) = \frac{C}{|x|^{N+2s}} + o(|x|^{-N-2s})$, as $|x| \to \infty$;
\item[(iv)] If $2 \leq p < 1 + (2s + \alpha)/N$ and $s > 1/2$, the Morse index of $u$ is
equal to one.
\end{itemize}
\end{proposition}
Take $p=2$ and $\alpha=2$ in \eqref{eq2:2022914}, by using the above Proposition, under condition $2<N<2+4s$,
we obtain a ground state $u \in H^s(\R^N)$ to $(\ref{eq:20220420-e1})$  which is
positive, regularity, radially symmetric and decreasing.

Now, we equivalent the equation $(\ref{eq:20220420-e1})$ as follows
  \beq\label{20221223-e0}
  \left\{
\begin{array}{ll}
\aligned
&(-\Delta)^s u+u=2uv~~~~\text{in}~~\R^N,\\
&-\Delta v=u^2~~~~\text{in}~~\R^N,\\
&(u, v)\in H,
\endaligned
\end{array}
\right. \quad\tag{Q}
\eeq
where $H:= H^s(\R^N)\times\dot{H}^1(\R^N).$
The ground state for \eqref{20221223-e0} can be obtained (up to scaling) by the constrained
minimization problem
\beq\label{eq:20220913-e1}
\aligned
\nu_s:=\frac{1}{3}&\inf\limits_{(u,v)\in H, \int_{\R^N}u^2 v dx=1}\left(\|u\|^2_{H^s(\R^N)}+\|v\|^2_{\dot{H}^1(\R^N)}\right)\\
=\frac{1}{3}&\inf\limits_{(u,v)\in H\setminus\{0\}}\frac{\|u\|^2_{H^s(\R^N)}+\|v\|^2_{\dot{H}^1(\R^N)}}{\int_{\R^N}u^2 v dx}.
\endaligned
\eeq
We point out that for any $(u,v)\in H$
$$\int_{\R^N} u^2 v dx\leq \left(\int_{\R^N} u^{\frac{4N}{N+2}} dx\right)^{\frac{N+2}{2N}}  \left(\int_{\R^N} v^{2^*}dx\right)^{\frac{1}{2^*}}<+\infty,$$
thanks to $\frac{4N}{N+2}\in[2, 2_s^*]$ and Sobolev embedding.

Let $(U_s, V_s) = (u_s(|x|), v_s(|x|))$ be a minimizer for $\nu_s$, i.e. $\int_{\R^N}U_s^2 V_s dx= 1$ and $3\nu_s = \|U_s\|_{H^s(\R^N)}^2+\|V_s\|^2_{\dot{H}^1(\R^N)}$, then it is a ground state solution of
  \beq\label{20221223-e11}
  \left\{
\begin{array}{ll}
\aligned
&(-\Delta)^s u+u=2\nu_s vu~~\hbox{in}~\R^N,\\
&-\Delta v=\nu_s u^2~~\hbox{in}~\R^N.
\endaligned
\end{array}
\right. \quad\tag{Q*}
\eeq
and so it solves $(\ref{eq:20220420-e1})$ (up to scaling).
Its derivatives $(\partial_i U_s, \partial_i V_s)$ are solution of the linearized
system
\beq\label{eq:20221118-e3}
  \left\{
\begin{array}{ll}
\aligned
&(-\Delta)^s\xi + \xi = 2\nu_s V_s\xi+2\nu_s \zeta U_s,\\
&-\Delta \zeta=2\nu_s \xi U_s.
\endaligned
\end{array}
\right.
\eeq
 Define the linear operator $\mathcal{L}_s^+$ associated to $(U_s, V_s)$ by
  \beq\label{eq:20220913-e2}
\mathcal{L}_s^+ \left(\begin{matrix} \xi\\ \zeta \end{matrix} \right)  =\left(\begin{matrix} (-\Delta)^s\xi + \xi- 2\nu_s V_s \xi - 2\nu_s\zeta U_s \\ -\Delta \zeta- 2\nu_s \xi U_s\end{matrix} \right).
\eeq
Now, we can state the nondegeneracy result.
\begin{theorem}\label{th1.2}
 Let $s\in(\frac{1}{4},1)$, $3\leq N<2+4s$. There exists $s_0\in(\frac{1}{4},1)$ such that for every $s\in(s_0,1)$, the operator $\mathcal{L}_s^+$ defined as in
\eqref{eq:20220913-e2} is nondegenerate around $(U_s, V_s)$. That is,
$$Ker(\mathcal{L}_s^+) = span\{(\partial_i U_s, \partial_i V_s), i=1,\cdot\cdot\cdot,N\}.$$
\end{theorem}

In the second part of this article, we prove the uniqueness of ground state.
\begin{theorem}\label{th1.3}
Let $s\in(\frac{1}{4},1)$, $3\leq N<2+4s$. There exists $s_0\in(\frac{1}{4},1)$ such that for every $s\in(s_0,1)$, the minimizer for
$\nu_s$ is unique, up to translations.
\end{theorem}

The paper is organized as follows. In Section 2 we provide uniform estimates and asymptotics of minimizers. In Section 3-4 we prove Theorem \ref{th1.2}-\ref{th1.3} respectively.

\begin{itemize}
\item[$\diamond$]
For the sake of simplicity, integrals over the whole $\R^N$ will be often written $\int$;
\item[$\diamond$] $\dot{H}^1(\R^N):=\left\{u\in L^{2^*}(\R^N): \nabla u\in L^2(\R^N)\right\}$, $\|u\|_{\dot{H}^1}=\|\nabla u\|_{L^2(\R^N)}$;
\item[$\diamond$] $\|\cdot\|_{H^s}$ denotes the standard norm for the fractional Sobolev space $H^s(\R^N)$;
\item[$\diamond$] $\|\cdot\|_{L^q}$ denotes the $L^{q}(\mathbb{R}^{N})$-norm for $q\in[1,\infty]$;
\item[$\diamond$] $o(1)$ denotes the infinitesimal as $n \to +\infty$;
\item[$\diamond$] $2^*_s:=\frac{2N}{N-2s}$ and $2^*:=\frac{2N}{N-2}$ denote fractional Sobolev critical exponent and Sobolev critical exponent, respectively;
\item[$\diamond$] $t'$ denotes the conjugate exponent of $t$, i.e., $\frac{1}{t'}+\frac{1}{t}=1$.
\item[$\diamond$] Unless otherwise specified, $C$ represents a pure constant independent of any variable.
\end{itemize}

\vskip4mm
{\section{ Uniform estimates and asymptotics } }
 \setcounter{equation}{0}

We call $\mathcal{M}_s$ the the space of these positive, radially symmetric minimizers $(u_s,v_s)\in H^s(\R^N)\times \dot{H}^1(\R^N)$ for $\nu_s$ normalized so that $\int u_s^2 v_s dx= 1$. Therefore if $(u_s,v_s)\in\mathcal{M}_s$ then
\beq\label{20221013-e0}
\|u_s\|_{L^\infty}= |u_s (x^s_0)|,\quad \|v_s\|_{L^\infty}= |v_s (x^s_0)|,
\eeq
for some $x^s_0\in\R^N$.
\begin{lemma}\label{3.1} $\nu_s$ is uniform bounded, i.e., $\sup\limits_{s\in(0,1]}
\nu_s < +\infty.$
\end{lemma}
\bp
Let $(u_1,v_1) \in \mathcal{M}_1.$ It follows from $|\xi|^{2s} \leq 1+|\xi|^2$ that
$$\|u_1\|^2_{H^s}+\|v_1\|^2_{\dot{H}^1} \leq 2\|u_1\|^2_{L^2}+\int_{\R^N} |\xi|^2 |\widehat{u}|^2 d\xi +\|v_1\|^2_{\dot{H}^1}\leq 2\|u_1\|^2_{H^1}+2\|v_1\|^2_{\dot{H}^1} =6\nu_1 . $$ Since $\nu _s \leq \|u_1\|^2_{H^s} +\|v_1\|^2_{\dot{H}^1} $, the desired result follows.
\ep

To get uniform bounds on the minimizers, we need the following regularity theory of fractional Laplacian.
\begin{proposition}\label{20221223-e2}(\cite[Proposition 2.1.9]{MR2270163}, see also \cite{MR3165278}.) Assume $s > 0$, $u\in L^\infty(\R^N)$ and $(-\Delta)^s u \in L^\infty(\R^N)$.
\\
 If $2s \leq 1$, then $u \in C^{0,\alpha}(\R^N)$ for any $\alpha < 2s$. Moreover
$$\|u\|_{ C^{0,\alpha}(\R^N )} \leq C(s, N, \alpha) \left(\|(-\Delta)^s u\|_{L^\infty} + \|u\|_{ L^\infty} \right).$$
 If $2s > 1$, then $u \in C^{1,\alpha}(\R^N)$ for any $\alpha < 2s - 1$. Moreover
$$\|u\|_{ C^{1,\alpha}(\R^N)} \leq C(s, N, \alpha) \left(\|(-\Delta)^s u\|_{L^\infty} + \|u\|_{ L^\infty} \right).$$
\end{proposition}

Now we use Gidas-Spruck's blowup method \cite{MR619749} to get uniform $L^\infty$ bounds on the minimizers.
\begin{lemma}\label{3.2} Given $s_0 \in (0, 1)$, we have
\beq\label{20221013-e1}
0 < \gamma_{s_0}:= \sup\limits_{s\in(s_0,1)}\sup\limits_{(u_s,v_s)\in\mathcal{M}_s}
(\|u_s\|_{L^\infty}+\|v_s\|_{L^\infty}) < \infty.
\eeq
Also, given $s_1 > 1/2$ and $\sigma\in(0, 1)$,
\beq\label{20221013-e11}
\sup\limits_{s\in(s_1,1)}\sup\limits_{(u_s,v_s)\in\mathcal{M}_s}\left(\|u_s\|_{C^{1,\sigma}(\R^N)}+\|v_s\|_{C^{1,\sigma}(\R^N)}\right) < \infty.
\eeq
\end{lemma}
\bp
$\gamma_{s_0}>0$ is obvious.
Now we prove the second inequality in \eqref{20221013-e1}. For this, we argue by contradiction: let
\beq\label{20221013-e2}
\lambda_s := \|u_s\|_{L^\infty}+\|v_s\|_{L^\infty}
\eeq
and assume that $\lambda_s \to+\infty$ for a sequence $s \to \bar{\sigma}\in [s_0, 1].$
We set
\beq\nonumber
w_s(x) := \lambda_s^{-1} u_s(\lambda_s^{\frac{2}{2s-N}}x+ x_0^s),
\quad \psi_s(x):=\lambda_s^{-1} v_s(\lambda_s^{\frac{2}{2s-N}}x+ x_0^s)
\eeq
so that
$$\|w_s\|_{L^\infty}+\|\psi_s\|_{L^\infty}=w_s(0)+\psi_s(0)=1,$$
and by $\frac{2s-2}{N-2s}<0$
$$\|w_s\|_{\dot{H}^s(\R^N)}+\|\psi_s\|_{\dot{H}^1(\R^N)}=\|u_s\|_{\dot{H}^s(\R^N)}+\lambda_s^{\frac{2s-2}{N-2s}}\|v_s\|_{\dot{H}^1(\R^N)}\leq \nu_s\leq C.$$
Therefore $(w_s,\psi_s)\rightharpoonup (w,\psi)$ in $H^t(\R^N )\times H^k(\R^N )$ for every $t <\bar{\sigma}$ and $k<1$, and
$$(w_s,\psi_s) \to (w,\psi)~\text{in}~L^2_{loc}(\R^N )\times L^2_{loc}(\R^N ),\quad \text{as}~s \to \bar{\sigma}.$$
Since $u_s$ satisfies \eqref{20221223-e11},
then $(w_s, \psi_s)$ satisfies
\beq\label{20221013-e3}
\left\{
\begin{array}{ll}
\aligned
&(-\Delta)^s w_s+ \lambda_s^{-\frac{4s}{N-2s}}w_s=2\nu_s\lambda_s^{1-\frac{4s}{N-2s}}\psi_s w_s\quad \text{in}~\R^N,\\
&-\Delta \psi_s =\nu_s\lambda_s^{1-\frac{4}{N-2s}} w_s^2\quad \text{in}~\R^N.
\endaligned
\end{array}
\right.
\eeq
By Proposition \ref{20221223-e2},
$$\|w_s\|_{C^{0,\alpha}(\R^N )} \leq C(s, N, \alpha) \left(\|(-\Delta)^s w_s\|_{L^\infty} + \|w_s\|_{L^\infty}\right)
,$$
$$\|\psi_s\|_{C^{0,\alpha}(\R^N )} \leq C(N, \alpha) \left(\|-\Delta \psi_s\|_{L^\infty} + \|\psi_s\|_{L^\infty}\right)
,$$
where one can fix $\alpha< 2 \bar{\sigma}$ for $2\bar{\sigma} < 1$ and $\alpha< 2\bar{\sigma}- 1$ for $2 \bar{\sigma}> 1$ and the constant
$C(s, N, \alpha)$ is bounded uniformly in $s \in [s_0, 1]$.
Then by \eqref{20221013-e3} and Lemma \ref{3.1}, we
see that $\|w_s\|_{C^{0,\alpha}(\R^N)}+\|\psi_s\|_{C^{0,\alpha}(\R^N)}$ is bounded uniformly when $s\to \bar{\sigma}$. Accordingly, by the Ascoli
theorem, we may suppose that $(w_s,\psi_s)$ converges locally uniformly to $(w,\psi)$ and passing to the
limit in \eqref{20221013-e3}, by $2s<N<2+4s$ we have that
\beq\nonumber
\left\{
\begin{array}{ll}
\aligned
&(-\Delta)^{\bar{\sigma}} w=0\quad \text{in}~\R^N,\\
&-\Delta \psi=0\quad \text{in}~\R^N.
\endaligned
\end{array}
\right.
\eeq
 Then by Liouville theorem we get $w=\psi \equiv 0$. In particular
$$0 = \lim\limits_{s\to\bar{\sigma}}(|w_s (0)|+|\psi_s(0)|) = \lim\limits_{s\to\bar{\sigma}}
\left(\lambda_s^{-1}|u_s (x^s_0)|+\lambda_s^{-1}|v_s (x^s_0)|\right)= 1,$$
due to \eqref{20221013-e0} and \eqref{20221013-e2}. This is a contradiction and so \eqref{20221013-e1} is proved.

 By using once again Proposition \ref{20221223-e2},
 for any $s\in (s_1, 1]$ and $\alpha<2s-1$,
$$\|u_s\|_{C^{1,\alpha}(\R^N )} \leq C(s, N, \alpha) \left(\|(-\Delta)^s u_s\|_{L^\infty} + \|u_s\|_{L^\infty}\right)
,$$
$$\|v_s\|_{C^{1,\alpha}(\R^N )} \leq C(N, \alpha) \left(\|-\Delta v_s\|_{L^\infty} + \|v_s\|_{L^\infty}\right)
,$$
where $C(s, N, \alpha)$ is uniformly bounded on $[s_1, 1]$. Then, \eqref{20221223-e11}, \eqref{20221013-e1} and Lemma \ref{3.1} imply \eqref{20221013-e11}.
\ep

\begin{corollary}\label{3.3} Given $s_0 \in (\frac{1}{4}, 1)$, we have
$$\sup\limits_{s\in(s_0,1)}\sup\limits_{(u_s,v_s)\in\mathcal{M}_s}
(\|u_s\|_{H^{2s}}+\|v_s\|_{\dot{H}^{2}}) < \infty.$$
\end{corollary}
\bp
First, from Lemma \ref{3.1} we have $$\|u_s\|^2_{H^s}+\|v_s\|^2_{\dot{H}^1} = 3\nu_s \leq C_1 ,$$ with $C_1>0$ independent of $s$ and $(u_s,v_s)$.

Let $s_0 \in (0,1)$,$(u_s,v_s)\in \mathcal{M}_s$ and
$$f_s^{(1)}(x):=2\nu _s v_s u _s(x)-u_s(x),\quad f_s^{(2)}(x):=\nu_s u^2_s(x).$$
By \eqref{20221013-e1}, Lemma \ref{3.1} and Lemma \ref{3.2}, we have
$$\int _{\R^N} \left|v_s u _s\right|^2dx \leq \|u_s\|_{L^2}^2\|v_s\|^2_{L^\infty}\leq C_2,$$
where $C_2>0$ is a constant independent of $s$ and $(u_s,v_s)$.  As a consequence,and using Lemma \ref{3.1} again, we obtain that
 $$\|f_s^{(1)}\|_{L^2}+\|f_s^{(2)}\|_{L^2} \leq 2\nu _s\|v_s u _s\|_{L^2}+\|u_s\|_{L^2}+\nu_s \|u_s\|_{L^2}^2\leq C_3,$$
 with $C_3>0$ independent of $s$ and $(u_s,v_s)$. Also, from \eqref{20221223-e11},
 \beq\nonumber
\left\{
\begin{array}{ll}
\aligned
&(-\Delta)^s u_s =f_s^{(1)},\\
&-\Delta v_s=f_s^{(2)},
\endaligned
\end{array}
\right.
\eeq
 that is, recalling \eqref{eq1:20221118},
 $$|\xi|^{2s} \widehat{u_s}=\widehat{f_s^{(1)}},\quad |\xi|^{2} \widehat{v_s}=\widehat{f_s^{(2)}}$$
 and so
 \beq\nonumber
 \aligned
 &\|u_s\|^2_{H^{2s}}+\|v_s\|^2_{\dot{H}^{2}} \\
 =&\|u_s\|^2_{L^2}+\int_{\R^N}|\xi|^{4s}|\widehat{u_s}|^2d\xi +\int_{\R^N}|\xi|^{4}|\widehat{v_s}|^2d\xi\\
 \leq &3\nu_s +\int_{\R^N}|\widehat{f_s^{(1)}}|^2 d\xi+\int_{\R^N}|\widehat{f_s^{(2)}}|^2 d\xi\\
 =&3\nu_s +\|f_s^{(1)}\|^2_{L^2}+\|f_s^{(2)}\|^2_{L^2}\leq C_1+C_3,
 \endaligned
 \eeq
 and the desired result follows.
\ep

\begin{proposition}\label{3.4} (\cite[Lemma 2.4]{FV}.) Let $s, \sigma \in (0, 1]$ and
\beq\label{eq:20221118-1e}
\delta > 2|\sigma - s|.
\eeq
For any $\xi \in \R^N \setminus \{0\}$,
\beq\label{eq:20221117-7}
||\xi|^{2s}-|\xi|^{2\sigma}|
\leq 4C_{\sigma,\delta}|\sigma-s|(1+|\xi|^{2(\sigma+\delta)}),
\eeq
and for any $\varphi\in H^{2(\sigma +\delta)}(\R^N )$,
$$\|(-\Delta)^\sigma\varphi - (-\Delta)^s\varphi\|_{L^2}
\leq C_{\sigma,\delta}|\sigma - s| \|\varphi\|_{H^{2(\sigma +\delta)}},$$
where
$$C_{\sigma,\delta}=\frac{1}{e}\left(\frac{1}{2\sigma}+\frac{1}{\delta}\right).$$
\end{proposition}

\begin{lemma}\label{3.5} Fix $\sigma\in(0, 1]$. Then $\lim\limits_{s\to\sigma}\nu_s = \nu_\sigma .$
\end{lemma}
\bp
Let $s_0 \in (0,1)$. Let $s$,$s' \in \left(s_0,1\right]$, that will be taken one close to the other, namely such that
\beq\nonumber
s>2|s-s'|.
\eeq
Let $(u_s, v_s) \in \mathcal{M}_s$. Since $\int u_s^2 v_s dx=1$, we get that $3\nu_{s'} \leq \|u_s\|^2_{H^{s'}}+\|v_s\|^2_{\dot{H}^{1}} $. Hence, by \eqref{eq:20221117-7}, we conclude that
\begin{align*}
3\nu_{s'}-3\nu_s &\leq \|u_s\|^2_{H^{s'}}-\|u_s\|^2_{H^s} \notag\\
&= \int_{\R^N}(|\xi|^{2s'}-|\xi|^{2s})|\widehat{u_s}|^2 \leq C(s_0)|s'-s|  \int_{\R^N}(1+|\xi|^{4s})|\widehat{u_s}|^2 \notag\\
&=C(s_0)|s'-s| \|u_s\|_{H^{2s}} \notag
\end{align*}
 Since the roles of $s$ and $s'$ may be interchanged, and recalling Corollary \ref{3.3}, we obtain that $$|\nu_{s'}-\nu_s| \leq C |s'-s|,$$ where the constant $C$ doesn't depend with $s', s$, and the desired result follows.
\ep

\begin{lemma}\label{20221103-l1} Fix $\sigma\in (0, 1]$. Let $s_n \in (\frac{1}{4}, 1)$ be such that $s_n \to\sigma$. Let $(u_{s_n}, v_{s_n})
\in\mathcal{M}_{s_n}$. Then
there exist $(\bar{u}, \bar{v})\in \mathcal{M}_\sigma$ and a subsequence (still denoted by $s_n$) such that
$$\|\psi_{s_n}\|_{\dot{H}^2}\to 0\quad \text{as}~~ n\to\infty;$$
 if $\sigma<1$ then
$$\|\omega_{s_n}\|_{H^{2s_n}}\to 0\quad \text{as}~~ n\to\infty;$$
 if $\sigma = 1$ then
$$\|\omega_{s_n}\|_{H^2}\to 0\quad \text{as}~~ n\to\infty,$$
where
\beq\label{eq:20221118_e2}
\omega_{s_n} (x) := u_{s_n} (x)-\bar{u},\quad \psi_{s_n} (x) := v_{s_n} (x)-\bar{v}.
\eeq
\end{lemma}

\bp
To alleviate the notation, we write $s$ instead of $s_n$. From Lemma \ref{3.2} and Corollary \ref{3.3} we have that $(u_s,v_s)$ is bounded in $H^t_{rad}(\R^N)\times H^k_{rad}(\R^N)$ for every $t<\sigma$ and $k<2$. Therefore, by Strauss's compactness embedding theorem \cite{MR454365}, we obtain that there exists $\bar{u}$ such that
$$u_s \rightarrow \bar{u} ~~\text{in}~L^q(\R^N)~for~every~q\in(2,2^{*}_{\sigma}),$$
$$v_s \rightarrow \bar{v}~~\text{in}~L^q(\R^N)~for~every~q\in(2,2^{*}).$$
Since we have uniform decay bounds at infinity and uniform $L^{\infty}$ bounds (recall Lemma \ref{3.2}), this and the interpolation inequality implies that the convergence also holds for $q\in \left(1,2\right]$, hence
\beq\label{eq:20221118-3}
\aligned
&u_s \rightarrow \bar{u} ~~\text{in}~L^q(\R^N)~for~every~q\in(1,2^{*}_{\sigma}),\\
&v_s \rightarrow \bar{v}~~\text{in}~L^q(\R^N)~for~every~q\in(1,2^{*}).
\endaligned
\eeq
In partical, $\bar{u}$ and $\bar{v}$ are radially symmetric. What is more, by Fatou's lemma, it follows that $\bar{v}\in \dot{H}^1(\R^N)$, $\bar{u} \in H^{\sigma}(\R^N)$ because $\int_{\R^N}|\xi|^{2s}|\widehat{u_s}|^2d \xi \leq 3\nu_s \leq C$. Also, by \eqref{20221223-e11},
\beq\label{eq:20221118-4}
\left\{
\begin{array}{ll}
\aligned
&\int_{\R^N} u_s (-\Delta)^s \varphi +\int_{\R^N} u_s \varphi =2\nu_s \int_{\R^N} v_s u _s \varphi~~~\forall \varphi \in C^{\infty}_c(\R^N),\\
&\int_{\R^N} u_s (-\Delta) \phi=\nu_s \int_{\R^N} u _s^2 \phi~~~\forall \phi \in C^{\infty}_c(\R^N).
\endaligned
\end{array}
\right.
\eeq
Using Proposition \ref{3.4},
$$\int_{\R^N} |\left[(-\Delta)^s\varphi-(-\Delta)^{\sigma}\varphi\right]|^2 dx \leq C(\sigma-s)^2 \int_{\R^N} (1+|\xi|^4)|\widehat{\varphi}|^2 d\xi \leq C\|\varphi\|^2_{H^2}.$$
Hence we can pass to the limit in \eqref{eq:20221118-4} and conclude that $(\bar{u}, \bar{v})$ is a weak solution to the system
\beq\label{eq:20221118-5}
\left\{
\begin{array}{ll}
\aligned
&(-\Delta)^{\sigma}\bar{u}+\bar{u}=2\nu_{\sigma} \bar{v} \bar{u}, \quad \text{in}~\R^N,\\
&-\Delta\bar{v}=\nu_\sigma \bar{u}^2,\quad \text{in}~\R^N,
\endaligned
\end{array}
\right.
\eeq
that belongs to $H^{\sigma}(\R^N)\times \dot{H}^{1}(\R^N)$.
So, by testing the equation against $(\bar{u}, \bar{v})$ itself, we see that
$$\|\bar{u}\|^2_{H^{\sigma}}+\|\bar{v}\|^2_{\dot{H}^{1}}=3\nu_{\sigma} \int \bar{u}^2\bar{v} =3\nu_{\sigma},$$
hence $(\bar{u},\bar{v})$ is a minimizer for $\nu_{\sigma}$.
Furthermore, by \eqref{20221223-e11}, \eqref{eq:20221118_e2} and \eqref{eq:20221118-5},
\beq\label{eq:20221118-6}
\left\{
\begin{array}{ll}
\aligned
&(-\Delta)^s \omega_s +\omega_s =
2\nu_s \left(\bar{v}\omega_s+\bar{u}\psi_s+\omega_s\psi_s\right)
+\left[(-\Delta)^{\sigma}\bar{u}-(-\Delta)^s\bar{u}\right]+2(\nu_s -\nu_{\sigma})\bar{v}\bar{u},\\
&-\Delta \psi_s =\nu_s (u_s+\bar{u})\omega_s.
\endaligned
\end{array}
\right.
\eeq
We observe that, since $\bar{u}, u_s, \bar{v}, v_s \in C^2(\R^N )$, \eqref{eq:20221118-6} holds pointwise and thus, from Proposition \ref{3.4} and \eqref{20221013-e1},
\beq\label{eq:20221224-8}
\aligned
&\|-\Delta \psi_s\|^2_{L^2} \leq \|\omega_s\|^2_{L^2}\rightarrow 0,~~\text{as}~s \nearrow \sigma,\\
&\|(-\Delta)^s \omega_s\|^2_{L^2} \leq \|\omega_s\|^2_{L^2}+C(|\sigma-s|^2+|\nu_{\sigma}-\nu_s|^2+\|\omega_s\|^2_{L^2}) \rightarrow 0,~~\text{as}~s \nearrow \sigma.
\endaligned
\eeq
\eqref{eq:20221224-8} implies that $\|\psi_s\|_{\dot{H}^{2}} \rightarrow 0$, $\|\omega_s\|_{H^{2s}} \rightarrow 0$ as $s \nearrow \sigma$, as desired.

Next we consider the case $\sigma=1$. By \eqref{20221223-e11} and Lemma \ref{3.2} we have that for every $s$ close to $1$ $$\|\partial_j u_s\|_{H^{2s}} \leq C.$$
From this, \eqref{eq:20221118-5} and  \eqref{eq:20221118_e2}, we deduce that
$$\|\partial_j \omega_s\|_{H^{2s}} \leq C.$$
In particular $\|\omega_s\|_{H^{2s+1}}$ is uniformly bounded. We let $f_s$ be the right hand side of \eqref{eq:20221118-6} so that
$$(-\Delta)^s \omega_s +\omega_s =f_s$$
and so
$$-\Delta\omega_s +\omega_s=f_s+\left[-\Delta\omega_s -(-\Delta)^s\omega_s\right].$$
Using Proposition \ref{3.4}, we conclude that, for every $\delta \in (0,\frac{1}{4})$ and $\delta>2|1-s|$,
$$\int_{\R^N} \left[(-\Delta)^1\omega_s -(-\Delta)^s\omega_s\right]^2 \leq C_{N,\delta}(1-s) \|\omega_s\|_{H^{2+2\delta}} \leq C(1-s)\|\omega_s\|_{H^{2s+1}} \leq C(1-s),$$
provided $s$ is close to $1$. Also, by recallig \eqref{eq:20221224-8} and \eqref{eq:20221118-3}, we obtain that $\|f_s\|_{L^2} \rightarrow 0$ as $s \nearrow 1$, and therefore $\|\omega_s\|_{H^{2}} \rightarrow 0$.
\ep

\vskip4mm
{\section{ The proof of Theorem \ref{th1.2}}}
 \setcounter{equation}{0}

{\subsection{ Local realization of $(-\Delta)^s$ for $s \in (0,1)$ }
Following \cite{CS,CGM,GZ}, we recall here an extension property that provides a local realization of the fractional Laplacian $(-\Delta)^s$ by means of a divergence operator in the upper half-space $\R_+^{N+1}$. Namely, given $u \in H^s(\R^N)$, there exists a unique $\mathcal{H}(u) \in H^1(\R^{N+1}_+;t^{1-2s})$ such that
\beq\label{eq:20221118-e1}
\left\{
\begin{array}{ll}
	\aligned
&div(t^{1-2s} \nabla \mathcal{H}(u))=0~~~\text{in}~\R^{N+1}_+,		 \\
&\mathcal{H}(u)=u~~~\text{in}~\R^N, \\
&\lim\limits_{t \searrow 0} t^{1-2s} \mathcal{H}(u)_t=d_s(-\Delta)^s u~~~\text{on}~\R^N,
\endaligned
\end{array}
\right.
\eeq
where $d_s =2^{2s-1}\frac{\Gamma(s)}{\Gamma(1-s)}$ is a positive normalization constant.
Equivalently for every $\Psi \in H^1(\R^{N+1}_+;t^{1-2s})$,
\beq\label{eq:20221118-e2}
\int_{\R^{N+1}_+}\nabla \mathcal{H}(u) \cdot \nabla \Psi t^{1-2s} dtdx=d_s \int_{\R^N} |\xi|^{2s}\widehat{u}\widehat{\Psi} d\xi,
\eeq
where here and hereafter we denote the trace of a function with the same letter. From now on, we use $\mathcal{H}$ to denote the $s$-harmonic operator. Moreover, for any $\Phi \in H^1(\R^{N+1}_+;t^{1-2s})$, the trace $\Phi$ on $\R^N$ belongs to $H^s(\R^N)$; and $\mathcal{H}(tr(\Phi)):=\mathcal{H}(\Phi)$ has minimal Dirichlet energy.

Hence $\mathcal{H}(u_s)$ is radially symmetric with respect to the $x$ variable and it is a minimizer for
\beq\label{eq:20221118-e3}
\aligned
\nu_s =\frac{1}{3}\inf\limits_{(U,V) \in (H^1(\R^{N+1}_+;t^{1-2s})\times \dot{H}^1(\R^{N}))\setminus\{0\}} \frac{d^{-1}_s \int_{\R^{N+1}_+} |\nabla U|^2 t^{1-2s}dtdx +\int_{\R^N} |U|^2
+\int_{\R^N} |\nabla V|^2}{\int_{\R^N}VU^2}.
\endaligned
\eeq
Moreover, by \eqref{eq:20221118-e1}, system \eqref{20221223-e11} become to
\beq\label{eq:20221118-e4}
\left\{
\begin{array}{ll}
	\aligned
	&div(t^{1-2s}\nabla\mathcal{H}(u))=0~~~\text{in}~\R^{N+1}_+ \\
	&d^{-1}_s t^{1-2s}\mathcal{H}(u)_t +\mathcal{H}(u) =2\nu_s v\mathcal{H}(u)~~~\text{on}~\R^N\\
&-div(\nabla v) =\nu_s\mathcal{H}(u)^2~~~\text{on}~\R^N.
\endaligned
\end{array}
\right.  	
\eeq

{\subsection{  Nondegeneracy }
In this section, we assume that $u_s\in\mathcal{M}_s$ and we prove
that it is nondegenerate for $s$ sufficiently close to $1$. For this, we denote by $\perp_{H^s}$ and $\perp_{\dot{H}^1}$ the
orthogonality relation in $H^s(\R^N)$ and $\dot{H}^1(\R^N)$ respectively.

From now on, we will use the uniqueness and nondegeneracy results for the Schr\"{o}dinger-Newton equation \beq\label{20221223-e111}
  \left\{
\begin{array}{ll}
\aligned
&-\Delta u+u=2\nu_1 vu~~\hbox{in}~\R^N,\\
&-\Delta v=\nu_1 u^2~~\hbox{in}~\R^N,
\endaligned
\end{array}
\right.
\eeq
where $N=3,4,5$, see \cite{WY,MR4204770}.
 Namely, we recall that there exists a unique radial minimizer $(U_1(x)=\bar{U}_1(|x|), V_1(x)=\bar{V}_1(|x|))$ for $\nu_1$, such that
\beq\label{eq:20221118-2}
Ker(\mathcal{L}_1^+)=span\{(\partial_jU_1,\partial_jV_1),j=1,\dots,N\}.
\eeq

Define the energy functional $I_s: H^s(\R^N)\times \dot{H}^1(\R^N) \mapsto \R$ for \eqref{20221223-e11} as
\beq\nonumber
I_s(u,v):=\frac{1}{2}\|u\|_{H^s}^2+\frac{1}{2}\|v\|_{\dot{H}^1}^2-\nu_s\int_{\R^N}u^2 v dx.
\eeq
Moreover,
$$ \langle I'_s(u,v),(\varphi,\psi)\rangle=(u,\varphi)_{H^s}+(v,\psi)_{\dot{H}^1}-2\nu_s\int_{\R^N}uv\varphi dx-\nu_s\int_{\R^N}u^2\psi dx.
$$
And the second order Gateaux derivative
$I''_s (u_s, v_s )$ possess the following property.
\bl\label{3.71} For every $\varphi \perp_{H^s} u_s$ and $\psi\perp_{\dot{H}^1} v_s$ we have that
\beq\nonumber
0 \leq I''_s (u_s, v_s )[(\varphi,\psi), (\varphi,\psi)] = \|\varphi\|_{H^s}^2+\|\psi\|_{\dot{H}^1}^2-2\nu_s\int v_s\varphi^2 dx-4\nu_s\int u_s\varphi \psi dx.
\eeq
\el
\bp
Let $\varepsilon>0$. Since $\varphi \perp_{H^s} u_s$ and $\psi\perp_{\dot{H}^1} v_s$, we have
\beq\label{eq:20221120-12}
\|\varepsilon \varphi+u_s\|^2_{H^s}={\varepsilon}^2 \|\varphi\|^2_{H^s} +\|u_s\|^2_{H^s}, \quad \|\varepsilon \psi+v_s\|^2_{\dot{H}^1}={\varepsilon}^2 \|\psi\|^2_{\dot{H}^1} +\|v_s\|^2_{\dot{H}^1}.
\eeq
Moreover, by using the system \eqref{20221223-e11}, we have
$$\int v_s u_s \varphi dx =0,\quad \int u_s^2 \psi dx=0.$$
Recalling also that functions in $\mathcal{M}_s$ are normalized with
$$\int v_s u_s^2 dx=1,$$ we obtain
\beq\label{eq:20221120-21}
\aligned
&\int(\varepsilon\psi+v_s) (\varepsilon\varphi+u_s)^2dx \\
=&1+\varepsilon^3\int \varphi^2\psi dx+
2\varepsilon^2\int\psi\varphi u_s dx+\varepsilon^2 \int\varphi^2 v_s dx+\varepsilon\int u_s^2\psi dx+2\varepsilon\int u_s\varphi v_s dx\\
=&1+\varepsilon^3\int \varphi^2\psi dx+
2\varepsilon^2\int\psi\varphi u_s dx+\varepsilon^2 \int\varphi^2 v_s dx.
\endaligned
\eeq
Now we recall the Taylor expansion
\beq\label{eq:20221120-41}
\frac{1}{1+x}=1-x +O(x^2)
\eeq
for small $x$. Thus, by inserting \eqref{eq:20221120-21} and \eqref{eq:20221120-41}, we obtain
$$\frac{1}{(\int(\varepsilon\psi+v_s)|\varepsilon \varphi+u_s|^{2})}=1-{\varepsilon}^2 \left(\int\varphi^2 v_s +2\int\psi\varphi u_s \right)+O({\varepsilon}^3).$$
From this and \eqref{eq:20221120-12} we obtain
\begin{align*}
&\frac{\|\varepsilon \varphi+u_s\|^2_{H^s}+\|\varepsilon \psi+v_s\|^2_{\dot{H}^1}}{\int(\varepsilon\psi+v_s)|\varepsilon \varphi+u_s|^{2}} \\
=&\left(1-{\varepsilon}^2 \left(\int\varphi^2 v_s +2\int\psi\varphi u_s\right)+O({\varepsilon}^3)\right)(\|u_s\|^2_{H^s}+\|v_s\|^2_{\dot{H}^1}+{\varepsilon}^2 \|\varphi\|^2_{H^s} +\varepsilon^2\|\psi\|^2_{\dot{H}^1} ) \notag\\
=&\|u_s\|^2_{H^s}+\|v_s\|^2_{\dot{H}^1} +{\varepsilon}^2\left( \|\varphi\|_{H^s}^2+\|\psi\|^2_{\dot{H}^1}-(\|u_s\|^2_{H^s}+\|v_s\|^2_{\dot{H}^1})\left(\int\varphi^2 v_s +2\int\psi\varphi u_s\right)\right)\\
&+O({\varepsilon}^3)\\
=&\|u_s\|^2_{H^s}+\|v_s\|^2_{\dot{H}^1} +{\varepsilon}^2\left( \|\varphi\|_{H^s}^2+\|\psi\|^2_{\dot{H}^1}-\nu_s\left(\int\varphi^2 v_s +2\int\psi\varphi u_s\right)
\right)+O({\varepsilon}^3).
\notag
\end{align*}
Then the desired result follows since $(u_s, v_s)$ attains the minimal value $\nu_s=\|u_s\|^2_{H^s}+\|v_s\|_{\dot{H}^1}^2$.
\ep

\bl\label{3.8} Let $\Phi\in H^1(\R^{N+1}_+; t^{1-2s} )$, $\Psi\in H^1(\R^{N} )$ be such that
\beq\label{eq:20221120-6}
\aligned
&d_s^{-1}\int_{\R^{N+1}_+}\nabla\Phi\cdot \nabla \mathcal{H}(u_s)t^{1-2s}dt dx +
\int_{\R^N}\Phi \mathcal{H}(u_s) dx = 0,\\
&\int_{\R^N}\nabla\Psi\cdot \nabla v_s dx = 0.
\endaligned
\eeq
Then
\beq
\aligned\label{eq:20221120-51}
&I_s''(\mathcal{H}(u_s),v_s)[(\Phi,\Psi),(\Phi,\Psi)] \\
=& d_s^{-1}\int_{\R^{N+1}_+}|\nabla\Phi|^2t^{1-2s}dz +\int_{\R^N}\Phi^2dx \\
&-2\nu_s\int_{\R^N}v_s\Phi^2dx-4\nu_s\int_{\R^N}\Phi \Psi u_sdx+\int_{\R^N} |\nabla\Psi|^2 dx \\
&\geq 0.
\endaligned
\eeq
In particular for any $(h,l)\in H^1(\R^2_{++};t^{1-2s}r^{N-1})\times \dot{H}^1(\R_+; r^{N-1})$
\beq\label{eq:20221120-4e}
\aligned
	&A_1((h,l),(h,l))\\
&:= \int_{\R^2_{++}} h^2_t t^{1-2s} r^{N-1}dtdr+(N-1)\int_{\R^2_{++}} h^2 t^{1-2s} r^{N-3}dtdr  +\int_{\R^2_{++}} h^2_r  t^{1-2s} r^{N-1}dtdr \\
&+d_s \int_{\R_+}h^2 r^{N-1}dr \\
	&-2\nu_s d_s \int_{\R_+}v_s h^2 r^{N-1}dr-4\nu_s d_s \int_{\R_+}h l u_s  r^{N-1}dr \\
&+ d_s\int_{\R_+} l_r^2 r^{N-1} dr+(N-1)d_s\int_{\R_+}l^2 r^{N-3}dr \\
&\geq 0.
\endaligned
\eeq
\el
\bp
The proof of \eqref{eq:20221120-51} is similar to the proof of Lemma \ref{3.71}, since $\mathcal{H}(u_s)$ minimizes \eqref{eq:20221118-e3}. Next, let $h\in H^1(\R^2_{++};t^{1-2s} r^{N-1})$, $l\in H^1(\R_+; r^{N-1})$ and define $\Phi_i(x):=h(t,|x|)\frac{x^i}{|x|}$, $\Psi_i(x):=l(|x|)\frac{x^i}{|x|}$.
Since
$$\partial_t\Phi_i=h_t\frac{x_i}{r},~~\partial_i\Phi_i=\frac{x_i}{r}\left(h_r\frac{x_i}{r}-h\frac{x_i}{r^2}\right)+h\frac{1}{r}, \quad \partial_j\Phi_i=\frac{x_i}{r}\left(h_r\frac{x_j}{r}-h\frac{x_j}{r^2}\right),~~j\neq i;$$
$$\partial_i\Psi_i=\frac{x_i}{r}\left(l_r\frac{x_i}{r}-l\frac{x_i}{r^2}\right)+l\frac{1}{r}, \quad \partial_j\Psi_i=\frac{x_i}{r}\left(l_r\frac{x_j}{r}-l\frac{x_j}{r^2}\right),~~j\neq i,$$
we have
\beq\label{eq:20221221-1}
\aligned
&\sum_{i} (\partial_t\Phi_i)^2=h_t^2; \\
&\sum_{i} (\partial_i\Phi_i)^2+\sum_{j\neq i}\sum_{i} (\partial_j\Phi_i)^2=h_r^2+(N-1)\frac{h^2}{r^2}; \\
&\sum_{i} (\partial_i\Psi_i)^2+\sum_{j\neq i}\sum_{i} (\partial_j\Psi_i)^2=l_r^2+(N-1)\frac{l^2}{r^2}.
\endaligned
\eeq
Since $\mathcal{H}(u_s)$ is radial in the $x$ variable, by odd symmetry we have
 $$\int_{\R^N}\nabla\Phi \nabla\mathcal{H}(u_s) dx=-\int_{\R^N}\Phi \Delta\mathcal{H}(u_s) dx=0 ,\quad \int_{\R^N}\Phi \mathcal{H}(u_s) dx=0 ,$$
 $$\int_{\R^N}\nabla\Psi \nabla v_s dx=-\int_{\R^N}\Psi \Delta v_s dx=0,$$
 and so $\Phi, \Psi$ satisfies \eqref{eq:20221120-6}. Then \eqref{eq:20221120-6},  \eqref{eq:20221120-51} and \eqref{eq:20221221-1} yield \eqref{eq:20221120-4e}.
\ep

By \cite{Lenzmann2009}, for linearized operators $L_+$ arising from ground states $Q$ for NLS with local nonlinearities,
it is a well-known fact that $Ker{L_+} = \{0\}$ when $L_+$ is restricted to radial functions implies that $Ker{L_+}$ is
spanned by $\{\partial_i Q\}_{i=1}^3$.

The proof, however, involves some Sturm-Liouville theory which is not applicable to $\mathcal{L}_+$ given in \eqref{eq:20220913-e2},
due to the presence of the nonlocal term. Also, recall that Newton's theorem \cite[(9.7.5)]{LiebLoss.2001} is not at our disposal,
since we do not restrict ourselves to radial functions anymore. To overcome this difficulty, we have
to develop Perron-Frobenius-type arguments for the action of $\mathcal{L}_+$ with respect to decomposition into
spherical harmonics.

Now we consider the spherical harmonics on $\R^N$ for $N \geq 2$, i.e., the solution of the classical eigenvalue problem
$$-\Delta_{\mathbb{S}^{N-1}}Y^i_k =\lambda_k Y^i_k~~~\text{on}~\mathbb{S}^{N-1}.$$
Let $n_k$ be the multiplicity of $\lambda_k$.
\begin{proposition}\label{20221223-e1}(\cite{GH}) The eigenvalue $\lambda_k=k(k + N - 2)$ for $k\in\mathbb{N}$.
$$n_0=1,~~ Y_0=Const;\quad \quad n_1=N, ~~Y^i_1=\frac{x^i}{|x|}~\text{for}~i=1,\ldots,N,$$
and
\begin{equation*}
\langle Y_k^i , Y_k^j \rangle_{L^2(\mathbb{S}^{N-1})} =
\left\{
\begin{array}{ll}
\aligned
&1, ~~~\text{if}~i = j;\\
&0, ~~~\text{if}~ i \neq j.
\endaligned
\end{array}
\right.
\end{equation*}
\end{proposition}

\bl\label{3.9} Let $(\varphi, \psi)\in \text{Ker}( I''_s (u_s, v_s ))$. Then
$$\varphi = \varphi_0(|x|) +\sum_{i=1}^N c^i\partial_i u_s,\quad
\psi = \psi_0(|x|) +\sum_{i=1}^N c^i\partial_i v_s, $$
where
$\varphi_0(r ) = \int_{\mathbb{S}^{N-1}}\varphi(r\theta)d\sigma(\theta)$, $\psi_0(r ) = \int_{\mathbb{S}^{N-1}}\psi(r\theta)d\sigma(\theta)$
and $c^i\in\R.$
\el
\bp
Let $(\varphi,\psi) \in Ker(I''_s(u_s,v_s))$ which means
\beq\label{20221219-e1}
    \left\{
\begin{array}{ll}
  \aligned
  &(-\Delta)^s \varphi +\varphi=2\nu_s \psi u_s +2\nu_s v_s \varphi,\\
  &-\Delta \psi=2\nu_s u_s \varphi.
  \endaligned
\end{array}
\right.
\eeq
 Let $\mathcal{H}(\varphi) \in H^1(\R^{N+1}_+;t^{1-2s})$ be the $s$-harmonic extension of $\varphi$ which satisfies
\beq\label{eq:20221120-7}
    \left\{
\begin{array}{ll}
\aligned
&d^{-1}_s \int_{\R^{N+1}_+} \nabla\mathcal{H}(\varphi) \cdot \nabla\Psi t^{1-2s}dtdx +\int_{\R^N} \varphi\Psi dx =2\nu_s\int_{\R^N}(\psi u_s+v_s\varphi)\Psi dx\\
&\int_{\R^{N}} \nabla\psi \cdot \nabla\Phi dx=2\nu_s\int_{\R^N} u_s \varphi\Phi dx
\endaligned
\end{array}
\right.
\eeq
for all $\Psi \in H^1(\R^{N+1}_+;t^{1-2s})$, $\Phi \in H^1(\R^{N})$. Now we decompose $\mathcal{H}(\varphi)$, $\psi$ in the spherical harmonics and we obtain
\beq\label{eq:20221120-8}
\mathcal{H}(\varphi)(t,x)=\sum \limits_{k \in \N} \sum \limits_{i=1}^{n_k} f^k_i(t,|x|)Y^i_k (\frac{x}{|x|}),\quad \psi(x)=\sum \limits_{k \in \N} \sum \limits_{i=1}^{n_k} g^k_i(|x|)Y^i_k (\frac{x}{|x|}),
\eeq
where $f^k_i \in H^1(\R^2_{++};t^{1-2s}r^{N-1})$, $g^k_i \in H^1(\R_+; r^{N-1})$.

\begin{remark}
Since $\varphi$ and $\psi$ are related, $f^k_i$ and $g^k_i$ are not independent.
\end{remark}

By testing the first equation in \eqref{eq:20221120-7} against the function $\Psi=h(t,|x|)Y^i_k$ and using polar coordinates and Proposition \ref{20221223-e1}, we obtain that, for any $h \in H^1(\R^2_{++}; t^{1-2s}r^{N-1})$, any $k \in\mathbb{N}$ and any $i \in \left[1,n_k\right]$,
\begin{equation}\label{20221219-ep1}
\aligned
&A_k((f^k_i,g^k_i), h)_1 \\
:=& \int_{\R^2_{++}} (f^k_i)_t h_t t^{1-2s} r^{N-1}dtdr +\int_{\R^2_{++}} (f^k_i)_r h_r t^{1-2s} r^{N-1}dtdr \\
&+\lambda_k \int_{\R^2_{++}} f^k_i ht^{1-2s} r^{N-3}dtdr +d_s \int_{\R_+} f^k_i hr^{N-1}dr \\
&-2\nu_s d_s \int_{\R_+} g^k_i u_s h  r^{N-1}dr-2\nu_s d_s \int_{\R_+} v_s f^k_i hr^{N-1}dr \\
=&0.
\endaligned
\end{equation}

By testing the second equation in  \eqref{eq:20221120-7} against the function $\Phi=l(|x|)Y^i_k$ and using polar coordinates and Proposition \ref{20221223-e1}, we obtain that, for any $l\in H^1(\R_+; r^{N-1})$, any $k \in\mathbb{N}$ and any $i \in \left[1,n_k\right]$,
\begin{equation}\label{20221221-ep1}
\aligned
&A_k((f^k_i,g^k_i),l)_2 \\
:=&\int_{\R_+} (g^k_i)_r l_r r^{N-1}dr+\lambda_k \int_{\R_+} g^k_i l r^{N-3}dr-2\nu_s \int_{\R_+} u_s f^k_i l r^{N-1}dr\\
=&0.
\endaligned
\end{equation}

Now we observe that
\begin{equation}\nonumber
\aligned
&A_k((f^k_i,g^k_i),(f^k_i,g^k_i)) \\
:=&A_k((f^k_i,g^k_i),f^k_i)_1+d_sA_k((f^k_i,g^k_i),g^k_i)_2 \\
=& \int_{\R^2_{++}} |(f^k_i)_t |^2 t^{1-2s} r^{N-1}dtdr +\int_{\R^2_{++}} |(f^k_i)_r |^2 t^{1-2s} r^{N-1}dtdr \\
&+\lambda_k \int_{\R^2_{++}} |f^k_i|^2t^{1-2s} r^{N-3}dtdr +d_s \int_{\R_+} |f^k_i|^2r^{N-1}dr \\
&-4\nu_s d_s \int_{\R_+}f^k_i g^k_i u_s   r^{N-1}dr-2\nu_s d_s \int_{\R_+} v_s |f^k_i|^2r^{N-1}dr\\
&+d_s\int_{\R_+} |(g^k_i)_r|^2 r^{N-1}dr+\lambda_k d_s\int_{\R_+} |g^k_i|^2 r^{N-3}dr \\
=&A_1((f^k_i,g^k_i),(f^k_i,g^k_i))+(\lambda_k -(N-1)) \int_{\R^2_{++}} |f^k_i|^2 t^{1-2s} r^{N-3}dtdr \\
&+(\lambda_k -(N-1))d_s \int_{\R_{+}} |g^k_i|^2  r^{N-3}dr\\
=&0.
\endaligned
\end{equation}
By Lemma \ref{3.8} ($A_1((f^k_i,g^k_i),(f^k_i,g^k_i))\geq0$) and the fact that $\lambda_k>N-1$ for $k \geq 2$, we obtain from the identities above that
\begin{align*}
0&=A_k((f^k_i,g^k_i),(f^k_i,g^k_i)) \notag \\
&\geq(\lambda_k -(N-1)) \int_{\R^2_{++}} |f^k_i|^2 t^{1-2s} r^{N-3}dtdr \\
&+(\lambda_k -(N-1)) d_s\int_{\R_{+}}|g^k_i|^2  r^{N-3}dr.\notag
\end{align*}
As a consequence, $f^k_i=0$ for every $k \geq 2$. Accordingly, \eqref{eq:20221120-8} becomes
$$\mathcal{H}(\varphi)(t,x)= \sum \limits_{i=1}^{N} f^1_i(t,|x|)Y^i_1 (\frac{x}{|x|}),\quad \psi(x)= \sum \limits_{i=1}^{N} g^1_i(|x|)Y^i_1 (\frac{x}{|x|}).$$
To complete the proof we need to characterize $f^1_i$ and $g^1_i$. For this, we notice that, for $i=1,\ldots,N$, the function
$$f^1_i(t,r)= \int_{\mathbb{S}^{N-1}} \mathcal{H}(\varphi)(t,r\theta){\theta}^i d\sigma(\theta),~~g^1_i(r)= \int_{\mathbb{S}^{N-1}} \psi(r\theta){\theta}^i d\sigma(\theta)$$ satisfies $f^1_i(t,0)=0$, $g^1_i(0)=0$,
\begin{equation}\label{eq:20221120-9}
\aligned
&A_1((f^1_i,g^1_i),h)_1 \\
=& \int_{\R^2_{++}} (f^1_i)_t h_t t^{1-2s} r^{N-1}dtdr +\int_{\R^2_{++}} (f^1_i)_r h_r t^{1-2s} r^{N-1}dtdr \\
&+(N-1) \int_{\R^2_{++}} f^1_i ht^{1-2s} r^{N-3}dtdr +d_s \int_{\R_+} f^1_i hr^{N-1}dr \\
&-2\nu_s d_s \int_{\R_+} g^1_i u_s h  r^{N-1}dr-2\nu_s d_s \int_{\R_+} v_s f^1_i hr^{N-1}dr \\
=&0,
\endaligned
\end{equation}
and
\begin{equation}\label{eq:20221120-92}
\aligned
&A_1((f^1_i,g^1_i),l)_2 \\
=&\int_{\R_+} (g^1_i)_r l_r r^{N-1}dr+(N-1)\int_{\R_+} g^1_i l r^{N-3}dr-2\nu_s\int_{\R_+} u_s f^1_i l r^{N-1}dr\\
=&0,
\endaligned
\end{equation}
for every $h \in H^1(\R^2_+;t^{1-2s}r^{N-1})$ and $l\in H^1(\R_+; r^{N-1})$, due to $\lambda_1=N-1$ and \eqref{20221219-ep1}-\eqref{20221221-ep1}.

Now we define $\bar{U}(t,|x|)=\mathcal{H}(u_s)(t,x)$ and $\bar{v}(|x|)=v_s(x)$.
Note that $\partial_r(r^{N-1}\partial_r \bar{U})=r^{N-1}\Delta_x \mathcal{H}(u_s)$ and
$$r^{N-1}{div}_{t,x}(t^{1-2s}\nabla_{t,x}\mathcal{H}(u_s))={div}_{t,r}(t^{1-2s}r^{N-1}\nabla_{t,r}\bar{U}).$$
Then we have
\beq\nonumber
\left\{
\begin{array}{ll}
	\aligned
    &{div}(t^{1-2s} r^{N-1} \nabla \bar{U})=0~~~\text{in}~\R^2_{++},	\\
	&\lim\limits_{t \searrow 0} -t^{1-2s} r^{N-1} \bar{U}_t+d_s r^{N-1} \bar{U}=d_s r^{N-1}2\nu_s \bar{U}v_s~~~\text{on}~\R_+, \\
&-\partial_{rr}\bar{v}-\frac{N-1}{r}\partial_r \bar{v}=\nu_s\bar{U}^2~~~\text{on}~\R_+, \\
	&\lim\limits_{r\searrow 0} r^{N-1} \bar{U}_r(t,0)=0,\quad \lim\limits_{r\searrow 0} r^{N-1} \bar{v}_r=0.
	\endaligned
\end{array}
\right.
\eeq
We differentiating the above equation with respect to $r$. We obtain
\beq\label{eq:20221120-10}
\left\{
\begin{array}{ll}
	\aligned
	&{div}(t^{1-2s} r^{N-1} \nabla \bar{U}_r)+(N-1)t^{1-2s}r^{N-3}\bar{U}_r=0~~~\text{in}~\R^2_{++},	\\
	&\lim\limits_{t \searrow 0} -t^{1-2s} r^{N-1} \bar{U}_{rt}+d_s r^{N-1} \bar{U}_r=2\nu_s d_s r^{N-1} (\bar{U}_r\bar{v}+\bar{U}\bar{v}_r)~~~\text{on}~\R_+, \\
&-\partial_r\left(\frac{1}{r^{N-1}}\partial_r(r^{N-1}\bar{v}_r)\right)=2\nu_s\bar{U}\bar{U}_r~~~\text{on}~\R_+, \\
	&\lim\limits_{r\searrow 0} r^{N-1} \bar{U}_r=0,\quad \lim\limits_{r\searrow 0} r^{N-1} \bar{v}_r=0.
	\endaligned
\end{array}
\right.
\eeq

By Proposition \ref{th1.1} (\cite[Theorem 1.1]{PGM}), $\bar{U}, \bar{v}$ are positive,
radially symmetric and decreasing, we may assume that $\bar{U}_r, \bar{v}_r<0$ on $\R^2_{++}$.

Given $f \in C^\infty_c (\R^2_{++}\bigcup\{t=0\})$, we define
$$\varsigma :=\frac{f}{\bar{U}_r} \in H^1(\R^2_{++};t^{1-2s}r^{N-1}).$$
Simple computations show that $$|\nabla f|^2=|V\nabla \varsigma|^2+\nabla \bar{U}_r \cdot \nabla(\bar{U}_r{\varsigma}^2).$$ Hence we have
\beq\label{20221220-ea1}
\aligned
&\int_{\R^2_{++}}|\nabla f|^2 t^{1-2s} r^{N-1}dtdr \\
=&\int_{\R^2_{++}}|\bar{U}_r\nabla\varsigma|^2 t^{1-2s} r^{N-1}dtdr +\int_{\R^2_{++}}\nabla (\bar{U}_r{\varsigma}^2) \cdot (t^{1-2s} r^{N-1} \nabla \bar{U}_r)dtdr.
\endaligned
\eeq
By testing the first equation of \eqref{eq:20221120-10} with $\frac{f^2}{\bar{U}_r}=\bar{U}_r{\varsigma}^2$, integrating by parts, and note that $\lim\limits_{t\to+\infty}t^{1-2s}\bar{U}_{rt}\frac{f^2}{\bar{U}_r}=\lim\limits_{r\to+\infty}r^{N-1}\bar{U}_{rr}\frac{f^2}{V}=\lim\limits_{r\to0^+}r^{N-1}\bar{U}_{rr}\frac{f^2}{\bar{U}_{r}}=0,$
we have
\beq\label{20221220-ea2}
\aligned
&(N-1)\int_{\R^2_{++}}f^2t^{1-2s}r^{N-3} dtdr \\
=&-\int_{\R^2_{++}}{div}(t^{1-2s}r^{N-1}\nabla \bar{U}_{r})\frac{f^2}{\bar{U}_{r}}dtdr \\
= &\lim\limits_{t\to0}\int_{\R_+}t^{1-2s}r^{N-1}\bar{U}_{rt} \frac{f^2}{\bar{U}_{r}}dr-\int_{\R^2_{++}}\nabla (\bar{U}_{r}{\varsigma}^2) \cdot (t^{1-2s} r^{N-1} \nabla \bar{U}_{r})dtdr.
\endaligned
\eeq
By testing the second equation of \eqref{eq:20221120-10} with $\frac{f^2}{\bar{U}_{r}}=\bar{U}_{r}{\varsigma}^2$,
\beq\label{20221220-ea3}
\aligned
&d_s\int_{\R_+} f^2r^{N-1} dr-d_s2\nu_s\int_{\R_+} (\bar{v}+\bar{U}\bar{v}_r/\bar{U}_{r})f^2r^{N-1} dr\\
=&-\lim\limits_{t\to0}\int_{\R_+}t^{1-2s}r^{N-1}\bar{U}_{rt} \frac{f^2}{\bar{U}_{r}}dr.
\endaligned
\eeq
By \eqref{20221220-ea1}-\eqref{20221220-ea3}, we get
\begin{equation}\label{202212201813}
\aligned
&\int_{\R^2_{++}}|\nabla f|^2 t^{1-2s} r^{N-1}dtdr+(N-1)\int_{\R^2_{++}}f^2 t^{1-2s} r^{N-3}dtdr\\
&+d_s\int_{\R_+}f^2 r^{N-1}dr
-d_s2\nu_s\int_{\R_+}(\bar{v}+\bar{U}\bar{v}_r/\bar{U}_{r}) f^2 r^{N-1}dr  \\
=&\int_{\R^2_{++}}|\bar{U}_{r}\nabla \varsigma|^2 t^{1-2s} r^{N-1}dtdr.
\endaligned
\end{equation}
Note from \eqref{eq:20221120-9} that
\begin{equation}\label{202212201812}
\aligned
&A_1((f,g),f)_1\\
=& \int_{\R^2_{++}} f^2 t^{1-2s} r^{N-1}dtdr +\int_{\R^2_{++}} f^2_r  t^{1-2s} r^{N-1}dtdr \\
&+(N-1) \int_{\R^2_{++}} f^2 t^{1-2s} r^{N-3}dtdr +d_s \int_{\R_+} f^2 r^{N-1}dr \\
&-2\nu_s d_s \int_{\R_+} fg \bar{U}  r^{N-1}dr-2\nu_s d_s \int_{\R_+} \bar{v} f^2 r^{N-1}dr.
\endaligned
\end{equation}
By \eqref{202212201813} and \eqref{202212201812}, we get
\begin{equation}\label{202212201814}
\aligned
&A_1((f,g),f)_1 \\
= &\int_{\R^2_{++}}|\bar{U}_{r}\nabla (f/ \bar{U}_r)|^2 t^{1-2s} r^{N-1}dtdr -2\nu_s d_s \int_{\R_+} fg \bar{U}  r^{N-1}dr \\
&+d_s 2\nu_s\int_{\R_+}\bar{U}\bar{v}_r/\bar{U}_{r} f^2 r^{N-1}dr.
\endaligned
\end{equation}

Given $g \in C^\infty_c (\R_{+})$,
by testing the third equation of \eqref{eq:20221120-10} with $\frac{g^2}{\bar{v}_{r}}r^{N-1}$, we have
\beq\nonumber
\aligned
&2\nu_s\int_{\R_+}\bar{U}\frac{\bar{U}_r}{\bar{v}_r} g^2 r^{N-1} dr \\
=&-\int_{\R_+}\partial_r\left(\frac{1}{r^{N-1}}\partial_r(r^{N-1}\bar{v}_r)\right)r^{N-1}\frac{g^2}{\bar{v}_r}dr \\
:=&I
\endaligned
\eeq
Integrating by parts, we get
\beq\nonumber
\aligned
I=&(N-1)\int_{\R_+}\frac{1}{r^{N-1}}\partial_r(r^{N-1}\bar{v}_r) r^{N-2}\frac{g^2}{\bar{v}_r}dr+\int_{\R_+} \frac{1}{r^{N-1}}\partial_r(r^{N-1}\bar{v}_r) r^{N-1} \partial_r \left(\frac{g^2}{\bar{v}_r}\right) dr \\
=&(N-1)\int_{\R_+} \partial_r(r^{N-1}\bar{v}_r)\frac{g^2}{r\bar{v}_r}dr+\int_{\R_+}\partial_r(r^{N-1}\bar{v}_r) \frac{2gg_r\bar{v}_r-g^2\bar{v}_{rr}}{\bar{v}_r^2} dr\\
=&-(N-1)\int_{\R_+} r^{N-1}\bar{v}_r \left(-\frac{1}{r^2} \frac{g^2}{\bar{v}_r}+\frac{1}{r} \frac{2gg_r\bar{v}_r-g^2\bar{v}_{rr}}{\bar{v}_r^2} \right) dr\\
&+\int_{\R_+}\left((N-1)r^{N-2}\bar{v}_r+r^{N-1}\bar{v}_{rr}\right)\frac{2gg_r\bar{v}_r-g^2\bar{v}_{rr}}{\bar{v}_r^2} dr\\
=&(N-1)\int_{\R_+} g^2 r^{N-3} dr-\int_{\R_+} \left[\left(\frac{\bar{v}_{rr}}{\bar{v}_r}g\right)^2-2\frac{\bar{v}_{rr}}{\bar{v}_r}gg_r\right] r^{N-1} dr .
\endaligned
\eeq
\allowbreak
Therefore,
\beq\label{202212201803}
\aligned
&2\nu_s\int_{\R_+}\bar{U}\frac{\bar{U}_r}{\bar{v}_r} g^2 r^{N-1} dr \\
=&(N-1)\int_{\R_+} g^2 r^{N-3} dr-\int_{\R_+} \left[\left(\frac{\bar{v}_{rr}}{\bar{v}_r}g\right)^2-2\frac{\bar{v}_{rr}}{\bar{v}_r}gg_r\right] r^{N-1} dr .
\endaligned
\eeq

Note from \eqref{eq:20221120-92} that
\begin{equation}\label{202212201804}
\aligned
A_1((f,g),g)_2=\int_{\R_+} g_r^2 r^{N-1}dr+(N-1)\int_{\R_+} g^2 r^{N-3}dr-2\nu_s \int_{\R_+} \bar{U} fg r^{N-1}dr.
\endaligned
\end{equation}
By \eqref{202212201803} and \eqref{202212201804}, we get
\begin{equation}\label{202212201805}
\aligned
&A_1((f,g), g)_2
=2\nu_s\int_{\R_+}\bar{U}\frac{\bar{U}_r}{\bar{v}_r} g^2 r^{N-1} dr-2\nu_s\int_{\R_+} \bar{U} fg r^{N-1}dr \\
&+\int_{\R_+} \left[g_r^2+\left(\frac{\bar{v}_{rr}}{\bar{v}_r}g\right)^2-2\frac{\bar{v}_{rr}}{\bar{v}_r}gg_r\right] r^{N-1} dr.
\endaligned
\end{equation}
Combining \eqref{202212201805} and \eqref{202212201814}, we get
\begin{equation}\label{202212201827}
\aligned
&A_1((f,g),(f,g))=A_1((f,g),f)_1+d_sA_1((f,g), g)_2 \\
=&\int_{\R^2_{++}}|\bar{U}_{r}\nabla (f/ \bar{U}_r)|^2 t^{1-2s} r^{N-1}dtdr +2\nu_s d_s \int_{\R_+}\left(\frac{\bar{v}_r}{\bar{U}_r}f^2+\frac{\bar{U}_r}{\bar{v}_r}g^2 -2fg  \right) \bar{U} r^{N-1}dr\\
&+d_s\int_{\R_+} \left[g_r^2+\left(\frac{\bar{v}_{rr}}{\bar{v}_r}g\right)^2-2\frac{\bar{v}_{rr}}{\bar{v}_r}gg_r\right] r^{N-1} dr .
\endaligned
\end{equation}
Therefore, by elementary inequality $$ a^2+b^2-2ab\geq0,$$
we obtain
\begin{equation}\label{202212201828}
\aligned
&A_1((f,g),(f,g)) \geq \int_{\R^2_{++}}|\bar{U}_{r}\nabla (f/ \bar{U}_r)|^2 t^{1-2s} r^{N-1}dtdr \\
+&2\nu_s d_s \int_{\R_+}\left(\sqrt{\frac{\bar{v}_r}{\bar{U}_r}}f-\sqrt{\frac{\bar{U}_r}{\bar{v}_r}}g\right)^2   \bar{U} r^{N-1}dr.
\endaligned
\end{equation}
In particular, by density we have that, for every $i=1,\ldots,N$,
\begin{equation}\nonumber
\aligned
0=&A_1((f^1_i,g^1_i),(f^1_i,g^1_i))\geq \int_{\R^2_{++}} \left|\bar{U}_{r}\nabla\left(\frac{f^1_i}{ \bar{U}_{r}}\right)\right|^2 t^{1-2s} r^{N-1}dtdr
\\
&+2\nu_s d_s \int_{\R_+}\left(\sqrt{\frac{\bar{v}_r}{\bar{U}_r}}f^1_i-\sqrt{\frac{\bar{U}_r}{\bar{v}_r}}g^1_i\right)^2   \bar{U} r^{N-1}dr.
\endaligned
\end{equation}
This implies that the last two terms vanish and therefore
$$\frac{f^1_i}{\bar{U}_{r}}=\frac{g^1_i}{\bar{v}_{r}} \equiv c^i$$
for some constant $c^i \in \R$. We then conclude that
$$f^1_i(0,|x|)=c^i \partial_r\bar{U}(0,|x|),\quad g^1_i(0,|x|)=c^i \partial_r\bar{v}(0,|x|) \quad\forall x \in \R^N.$$
Thus, we have proved that for any $(\varphi,\psi) \in Ker(I''_s(u_s,v_s))$
$$\mathcal{H}(\varphi)(0,x)=\varphi(x)=f^0_1(0,|x|)+\sum \limits_{i=1}^{N} f^1_i(0,|x|) \frac{x^i}{|x|}=f^0_1(0,|x|)+\sum \limits_{i=1}^{N} c^i \partial_i u_s(x),$$
and
$$\psi(x)=g^0_1(|x|)+\sum \limits_{i=1}^{N} g^1_i(|x|) \frac{x^i}{|x|}=g^0_1(|x|)+\sum \limits_{i=1}^{N} c^i \partial_i v_s(x),$$
as desired.

\ep

Now we are ready to prove our nondegeneracy result for $s$ close to $1$.

{\subsection{  Completion of the proof of Theorem \ref{th1.2}. }
 Let $(w_s, \vartheta_s)  \in Ker(I''_s(u_s, v_s ))$ be a radial
function. To proof Theorem \ref{th1.2}, according to Lemma \ref{3.9}, it is suffice to prove the following Claim.

Claim: If $s$ is close to 1, we have $w_s =\vartheta_s \equiv 0$.

Assume by contradiction that there exists a sequence $s_n$-still denoted by $s$- with $s\nearrow 1$ and such that $(w_s, \vartheta_s)\neq (0,0)$. Up to normalization, we can assum that $\int (I_2\star w_s^2) w_s^2 dx=1$. By Lemma \ref{3.1} and Lemma \ref{3.2}, using equation \eqref{20221219-e1} we have
\beq\label{20230404-e1}
\aligned
&\|w_s\|^2_{H^s} = 2\nu_s\int v_s w_s^2dx+2\nu_s\int u_sw_s\vartheta_sdx =(2\nu_s)^2\int (I_2\star u_s^2)w_s^2dx+(2\nu_s)^2\int (I_2\star u_sw_s)u_sw_sdx ; \\
&\|\vartheta_s\|^2_{\dot{H}^1}=2\nu_s\int u_sw_s\vartheta_s dx=(2\nu_s)^2\int (I_2\star u_sw_s)u_sw_sdx.
\endaligned
\eeq
In order to prove the boundedness of $w_s$ in $H^s(\R^N)$, we need some properties of the Riesz potential $I_\alpha$. First, by the semigroup property of the Riesz potential, we have the following conclusion.
\begin{proposition}\label{20221124-p1} (\cite[Proposition 2.2]{luo2022bifurcation}) Let $\alpha\in(0,N)$ and $f, g \in L^{\frac{2N}{N+\alpha}}(\R^N)$. Then
\beq\int
(I_\alpha\star f)gdx =
\int(I_{\frac{\alpha}{2}}\star f)
(I_{\frac{\alpha}{2}}\star g)dx=\int
(I_\alpha\star g)f
dx.
\eeq
Moreover,
\beq\label{20221124-e1}
\int(I_\alpha\star f)gdx \leq
\left[\int(I_\alpha\star f)fdx\right]^{\frac{1}{2}}\left[\int(I_\alpha\star g)gdx\right]^{\frac{1}{2}}.
\eeq
\end{proposition}
\begin{corollary}\label{20221124-c1} Let $u, v \in L^{\frac{4N}{N+2}}(\R^N)$.
$$\int(I_2\star uv)uv dx\leq \int(I_2\star u^2 )v^2dx \leq \left(\int(I_2\star u^2 )u^2dx\right)^{\frac{1}{2}}\left(\int(I_2\star v^2 )v^2dx\right)^{\frac{1}{2}}.$$
\end{corollary}
\bp
For any $x\in\R^N$, by the H\"{o}lder inequality, one has
$$\left|I_2\star uv \right|(x) \leq (I_2\star u^2 )^{\frac{1}{2}}(I_2\star v^2 )^{\frac{1}{2}}(x).$$
Then, using the H\"{o}lder inequality again,  one has
$$\int(I_2\star uv)uv dx \leq \left(\int(I_2\star u^2 )v^2dx\right)^{\frac{1}{2}}\left(\int(I_2\star v^2 )u^2dx\right)^{\frac{1}{2}}.$$
Therefore, by Proposition \ref{20221124-p1}, we get the conclusion.
\ep

Using \eqref{20230404-e1}, Proposition \ref{20221124-p1}, Corollary \ref{20221124-c1} and $\int v_su_s^2dx=\int(I_2\star u_s^2)u_s^2dx=1$, one has
\beq\label{eq:20221121-1}
\|w_s\|^2_{H^s}+ \|\vartheta_s\|^2_{\dot{H}^1}\leq Const.
\eeq
By \eqref{eq:20221121-1}, $w_s$ is a radial sequence and bounded in $H^t(\R^N)$ for every $t \in (\frac{1}{2},1)$. Then by Strauss's compactness embedding, up to a subsequence,
$$w_s \rightarrow w~~\text{in}~~L^q(\R^N)~\forall q \in (2,2^*),\quad \text{as}~s\to1.$$
Recalling Lemma \ref{3.1} ($\|u_s\|_{H^s}\leq\nu_s\leq C$), we have (up to a subsequence) $u_s \rightarrow U_1$ in $L^q(\R^N)$ for every $q \in (2,2^*)$.
Since we have uniform decay bounds at infinity and uniform $L^{\infty}$ bounds (recall Lemma \ref{3.2}), this and the interpolation inequality implies that the convergence also holds for $q\in \left(1,2\right]$. That is
$$ u_s \rightarrow U_1~~\text{in}~~ L^q(\R^N)~ \forall q \in (1,2^*),\quad \text{as}~s\to1. $$
On the other hand, since $\vartheta_s, v_s$ are bounded in $\dot{H}^1(\R^N) \subset H^1_{loc}(\R^N)$, we have (up to a subsequence)
$$\vartheta_s \to \vartheta,~~v_s \to V_1~~\text{in}~~L^q_{loc}(\R^N)~ \forall q \in (1,2^*),\quad \text{as}~s\to1.$$
Then by the H\"{o}lder inequality and $2\in (1,2^*)$, for any $\varphi \in C^\infty_c (\R^N)$ we have
\beq\nonumber
\aligned
 \lim\limits_{s\to1}\int_{\R^N}u_s\vartheta_s \varphi =\int_{\R^N}U_1\vartheta \varphi.
\endaligned
\eeq
By using the H\"{o}lder inequality with conjugate exponentials: $2+\epsilon\in (2,2^*)$ and $(2+\epsilon)'\in(1,2)$, for any $\varphi, \phi \in C^\infty_c (\R^N)$ we have
\beq\nonumber
\aligned
 \lim\limits_{s\to1}\int_{\R^N} w_sv_s\varphi=\int_{\R^N}  w V_1\varphi,\quad
 \lim\limits_{s\to1}\int_{\R^N} w_s u_s\phi=\int_{\R^N}  w U_1\phi.
\endaligned
\eeq
In particular, by $\dot{H}^1(\R^N)\hookrightarrow L^{2^*}(\R^N)$ we have $\vartheta_s \to \vartheta$ in $L^{2^*}(\R^N)$. By $\frac{4N}{N+2}\in(1, 2^*)$ we have $w_s^2\to w^2$ in $L^{\frac{2N}{N+2}}(\R^N)$.
 Next we observe that $(w_s,\vartheta_s)$ is a solution of the linearized equation and therefore for any $\varphi, \phi \in C^\infty_c (\R^N)$
\beq\nonumber
  \left\{
\begin{array}{ll}
\aligned
&\int_{\R^N} w_s (-\Delta)^s \varphi dx+\int_{\R^N} w_s \varphi dx=2\nu_s\int_{\R^N}w_sv_s\varphi dx+2\nu_s\int_{\R^N}u_s\vartheta_s \varphi dx, \\
&-\int_{\R^N} \vartheta_s \Delta \phi dx=2\nu_s\int_{\R^N}w_s u_s \phi dx,
\endaligned
\end{array}
\right.
\eeq
so by \eqref{20221013-e1}, Lemma \ref{3.5}, $w_s \rightarrow w$ in $L^2_{loc}(\R^N)$ and the fact that $(-\Delta)^s \varphi \rightarrow -\Delta \varphi $ in $L^2(\R^N)$ thanks to Proposition \ref{3.4}, we infer that
\beq\nonumber
  \left\{
\begin{array}{ll}
\aligned
&-\int_{\R^N} w \Delta \varphi dx+\int_{\R^N} w\varphi dx=2\nu_1\int_{\R^N}w V_1\varphi dx+2\nu_1\int_{\R^N}U_1 \vartheta \varphi dx, \\
&-\int_{\R^N} \vartheta \Delta \phi dx=2\nu_1\int_{\R^N}w U_1 \phi dx,
\endaligned
\end{array}
\right.
\eeq
 Applying Fatou lemma to \eqref{eq:20221121-1}, we get $w\in H^1(\R^N)$ and $\vartheta\in \dot{H}^1(\R^N)$. We then conclude that $(w,\vartheta)$ is radial, nontrivial and belongs to $Ker(I''_1(U_1,V_1))$. This is clearly a contradiction and the claim is proved.

\vskip4mm
{\section{ The proof of Theorem \ref{th1.3}}}
\setcounter{equation}{0}

{\subsection{ Convert the system \eqref{20221223-e11} into a single equation }

The system \eqref{20221223-e11} equivalent to
\beq\label{eq:20220913-e3}
(-\Delta)^su+ u = 2\nu_s^2 (I_2\star u^2)u.\quad\tag{P*}
\eeq
The ground state for \eqref{eq:20220913-e3} equivalent to the minimizer for
\beq\label{20221224-i1}
2\nu_s^2=\inf\limits_{u\in H^s(\R^N)\setminus\{0\}}\frac{\|u\|_{H^s}^2}{\|u\|_{HL}^2},
\eeq
up to scaling, where
$$\|u\|_{HL}:=\left(\int(I_2\star u^2)u^2dx\right)^{\frac{1}{4}}.$$
Let $\mathcal{M}_s$ be still defined as the space of these positive, radially symmetric minimizers $u_s\in H^s(\R^N)$ for
\eqref{20221224-i1} normalized so that $\|u_s\|_{HL}= 1$.

The energy functional $J_s$ associated to \eqref{eq:20220913-e3} defined as
$$J_s(u,\nu):=\frac{1}{2}\|u\|_{H^s}^2-\frac{\nu^2}{2}\|u\|_{HL}^4.$$
To study $J_s$ in $H^s(\R^N)$, we need the Hardy--Littlewood--Sobolev inequality (or abbreviated H-L-S inequality).
\begin{proposition}\cite{MR3625092}  Let $t$, $r>1$ and $0<\alpha<N$ with $\frac{1}{t}+\frac{1}{r}=1+\frac{\alpha}{N}$, $f\in L^{t}(\R^N)$ and $h\in L^{r}(\R^N)$. There exists a constant $C(N,\alpha,t,r)$, independent of $f,h$, such that
\begin{equation*}
\|I_\alpha\star h\|_{L^{t'}}\leq C(N,\alpha,t,r)\|h\|_{L^r}
\end{equation*}
and
\begin{equation*}
\int\left(I_\alpha\star h\right)f dx\leq C(N,\alpha,t,r) \|f\|_{L^t}\|h\|_{L^r}.
\end{equation*}
\end{proposition}
Combining the H-L-S inequality and the Sobolev inequality yields
\beq\label{20221124-e2}
\|u\|_{HL}^4 \leq C|u|^4_{L^{\frac{4N}{N+2}}}\leq C\|u\|^4_{H^s},\quad N<4s+2\Rightarrow \frac{4N}{N+2}<\frac{2N}{N-2s}.
\eeq
Furthermore, in \cite{MMV} the authors noted that $\|u\|_{HL}^4$ is naturally settled in the so called Coulomb
spaces $L^{HL}$, defined as the vector spaces of measurable functions $u : \R^N \to \R$ such that $\|u\|_{HL}$ is
finite. They also proved that the quantity $\|u\|_{HL}$
defines a norm, which will guarantees the convexity of the functional $\|\cdot\|_{HL}^4$. Hence, inequality
\eqref{20221124-e2} corresponds to the embedding $$H^s(\R^N)\subset L^{\frac{4N}{N+2}}(\R^N) \subset L^{HL}.$$ The paper \cite{MMV} then introduces and
carefully studies the Couloumb-Sobolev spaces and regularity properties in this framework.

{\subsection{ Preliminary observations }

First, we observe that for every $\varphi\in H^s(\R^N)$,
\beq\label{eq:20221121-4}
J''_s (u_s, \nu_s )[\varphi, \varphi] = \|\varphi\|_{H^s}^2-2\nu_s^2\int(I_2\star u_s^2)\varphi^2 dx-4\nu_s^2\int(I_2\star u_s\varphi)u_s\varphi dx.
\eeq
Similar to $I''_s (u_s, v_s )$, $J''_s (u_s, \nu_s )$ also have the following property.
\bl\label{3.7} For every $\varphi \perp_{H^s} u_s$,
$$ J''_s (u_s, \nu_s )[\varphi, \varphi]\geq 0.$$
\el
\bp
Let $\varepsilon>0$. Since $\varphi \perp_{H^s} u_s$, we have
\beq\label{eq:20221120-1}
\|\varepsilon \varphi+u_s\|^2_{H^s}={\varepsilon}^2 \|\varphi\|^2_{H^s} +\|u_s\|^2_{H^s}.
\eeq
We observe that
\beq\label{eq:20221120-2}
\aligned
&\int(I_2\star|\varepsilon \varphi+u_s|^{2})|\varepsilon \varphi+u_s|^{2}dx \\
=&\int(I_2\star|u_s|^{2})|u_s|^{2}dx+\varepsilon^4\int(I_2\star|\varphi|^{2})|\varphi|^{2}dx +2\varepsilon^2\int(I_2\star|u_s|^{2})|\varphi|^{2}dx \\
&+4\varepsilon\int\left(I_2\star u_s^2\right) u_s \varphi dx +4\varepsilon^3\int\left(I_2\star \varphi^2\right) u_s \varphi dx+4\varepsilon^2\int\left(I_2\star u_s\varphi\right) u_s \varphi dx.
\endaligned
\eeq
Furthermore, by testing \eqref{eq:20220913-e3} against $\varphi$ and using again that $\varphi \perp_{H^s} u_s$, we conclude that $$\int\left(I_2\star u_s^2\right) u_s \varphi dx=0,$$ hence the first order in $\varepsilon$ in \eqref{eq:20221120-2} vanishes. Consequently, recalling also that functions in $\mathcal{M}_s$ are normalized with $\|u\|_{HL}=1$, we write \eqref{eq:20221120-2} as
\beq\label{eq:20221120-3}
\aligned
&\int\left(I_2\star |\varepsilon \varphi+u_s|^{2}\right)|\varepsilon \varphi+u_s|^{2}dx \\
=&1+2\varepsilon^2\left(2\int\left(I_2\star u_s\varphi\right) u_s \varphi dx+\int\left(I_2\star u_s^2\right) \varphi^2\right)dx +O({\varepsilon}^3).
\endaligned
\eeq
Now we recall the Taylor expansion
\beq\label{eq:20221120-4}
\frac{1}{(1+x)^{\frac{1}{2}}}=1-\frac{1}{2}x +O(x^2)
\eeq
for small $x$. Thus, by inserting \eqref{eq:20221120-3} and \eqref{eq:20221120-4}, we obtain
\beq\nonumber
\aligned
&\frac{1}{\left(\int\left(I_2\star |\varepsilon \varphi+u_s|^{2}\right)|\varepsilon \varphi+u_s|^{2}dx\right)^{\frac{1}{2}}} \\
=&1-{\varepsilon}^2 \left(2\int\left(I_2\star u_s\varphi\right) u_s \varphi dx+\int\left(I_2\star u_s^2\right) \varphi^2 dx\right)+O({\varepsilon}^3).
\endaligned
\eeq
 From this and \eqref{eq:20221120-1} we obtain
\begin{align*}
&\frac{\|\varepsilon \varphi+u_s\|^2_{H^s}}{(\int\left(I_2\star |\varepsilon \varphi+u_s|^{2}\right)|\varepsilon \varphi+u_s|^{2}dx)^{\frac{1}{2}}} \\
&=\left(1-{\varepsilon}^2 \left(2\int\left(I_2\star u_s\varphi\right) u_s \varphi dx+\int\left(I_2\star u_s^2\right) \varphi^2 dx\right)+O({\varepsilon}^3)\right)({\varepsilon}^2 \|\varphi\|^2_{H^s} +\|u_s\|^2_{H^s}) \notag\\ &=\|u_s\|^2_{H^s} +{\varepsilon}^2\left( \|\varphi\|_{H^s}^2-\|u_s\|^2_{H^s}\int(I_2\star u_s^2)\varphi^2 dx-2\|u_s\|^2_{H^s}\int(I_2\star u_s\varphi)u_s\varphi dx
\right)+O({\varepsilon}^3). \notag
\end{align*}
Since $u_s$ attains the minimal value $2\nu_s^2=\|u_s\|^2_{H^s}$, thus $ J''_s (u_s, \nu_s )[\varphi, \varphi]\geq 0.$
\ep

\bl\label{4.1} Let $\Lambda_s := (Ker(J''_s (u_s, \nu_s )) \oplus \R u_s )^{\bot_{H^s}}.$
\\
(i) We have
\beq\label{eq:20221121-2}
J''_s (u_s, \nu_s )[u_s, u_s] = -2\|u_s\|^2_{H^s}.
\eeq
(ii) There exists $s_0 \in (0, 1)$ such that for every $s \in (s_0, 1)$ and every minimizer $u_s$ for $2\nu_s^2$
\beq\label{eq:20221121-3}
K(s, u_s ) := \inf\limits_{\varphi\in\Lambda_s\setminus\{0\}}
\frac{J''_s (u_s, \nu_s )[\varphi, \varphi]}{\|\varphi\|_{H^s}^2}> 0.
\eeq
(iii) Let
$$\Lambda^r_s:= \left\{\varphi\in H^s_{rad}(\R^N), \quad \varphi \bot_{H^s} u_s \right\}$$
and
$$K_r(s, u_s ) := \inf\limits_{\varphi\in\Lambda_s^r\setminus\{0\}}
\frac{J''_s (u_s, \nu_s )[\varphi,\varphi]}{\|\varphi\|_{H^s}^2}.$$
Then there exits $s_0 \in (0, 1)$ such that
\beq\label{eq:20221122-1}
\inf\limits_{s\in(s_0,1]}\inf\limits_{u\in\mathcal{M}_s}
K_r (s, u) > 0.
\eeq
\el
\bp
(i). Using \eqref{eq:20221121-4}, we immediately get \eqref{eq:20221121-2}.\par

(ii). We first show that for any $\varphi \in \Lambda_s$
\beq\label{eq:20221121-5}
J''_s(u_s,\nu_s)[\varphi,\varphi]=0 \Rightarrow \varphi \equiv 0.
\eeq
That is to say that $J''_s(u_s,\nu_s)$ defines a scalar product on $\Lambda_s$ by Lemma \ref{3.7}. For this, assume that $\varphi \in \Lambda_s$ and
$$J''_s(u_s,\nu_s)[\varphi,\varphi]=0.$$
Pick $\psi \in H^s(\R^N)$ such that $\psi \perp_{H^s} u_s$. We also have $\varphi \perp_{H^s} u_s$ due to $\varphi \in \Lambda_s$. Then by Lemma \ref{3.7}, we have
$$J''_s(u_s,\nu_s)[\varphi +\varepsilon \psi,\varphi +\varepsilon \psi] \geq 0.$$
Hence for any $\varepsilon\in\R$
\begin{align*}
0&\leq J''_s(u_s,\nu_s)[\varphi,\varphi]+2\varepsilon J''_s(u_s,\nu_s)[\varphi,\psi]+ {\varepsilon}^2 J''_s(u_s,\nu_s)[\psi,\psi] \notag \\
&=2\varepsilon J''_s(u_s,\nu_s)[\varphi,\psi] +{\varepsilon}^2 J''_s(u_s,\nu_s)[\psi,\psi].\notag
\end{align*}
Then we conclude that
\beq\label{eq:20221121-6}
J''_s(u_s,\nu_s)[\varphi,\psi]=0~~~\text{for any}~\psi \perp_{H^s} u_s.
\eeq
Now we observe from $\varphi \perp_{H^s} u_s$ and \eqref{eq:20220913-e3} that
$$0=\left\langle \varphi, u_s \right\rangle_{H^s} =2\nu_s^2 \int_{\R^N} \left(I_2\star u_s^2\right)u_s \varphi dx,$$
and so
$$J''_s(u_s,\nu_s)[\varphi,u_s]=\left\langle \varphi, u_s \right\rangle_{H^s} -6\nu_s^2 \int_{\R^N}\left(I_2\star u_s^2\right) u_s \varphi dx=0.$$
This and \eqref{eq:20221121-6} yield $\varphi \in Ker(J''_s(u_s,\nu_s))$. Since also $\varphi \perp_{H^s} Ker(J''_s(u_s,\nu_s))$ it follows that $\varphi =0$, and \eqref{eq:20221121-5} is proved.

Now we prove \eqref{eq:20221121-3} by contradiction. Assume that exits a sequence $\varphi_n \in \Lambda_s$ such that $\|\varphi_n\|_{H^s}=1$ and
\beq\label{eq:20221121-7}
J''_s(u_s,\nu_s)[\varphi_n,\varphi_n] \rightarrow 0~~~\text{as}~n\rightarrow\infty.
\eeq
Let $\varphi$ be the weak limit of $\varphi_n$ in $H^s(\R^N)$. Then, by Lemma \ref{3.7}, we have that
$$0\leq J''_s(u_s,\nu_s)[\varphi,\varphi] \leq \liminf J''_s(u_s,\nu_s)[\varphi_n,\varphi_n]=0.$$
We deduce from this and \eqref{eq:20221121-5} that $\varphi=0$, that is
\beq\label{eq:20221121-8}
\varphi_n\rightharpoonup0~\text{ in}~H^s(\R^N).
\eeq
Now, since $u_s \in L^{HL}$, given $\varepsilon \geq 0$ there exists $w_\varepsilon \in C^\infty_c (\R^N)$ such that
\beq\label{eq:20221121-9}
\|u_s -w_\varepsilon\|_{HL} < \varepsilon
\eeq
By \eqref{eq:20221121-8} and the compact embedding in fractional Sobolev spaces ( \cite[Theorem 7.1]{DPV}), we obtain
$$\varphi_n\to 0~\text{ in}~L^2_{loc}(\R^N)$$
 and therefore
\beq\label{eq:20221121-10}
\left|\int_{\R^N} \left(I_2\star w_\varepsilon^2\right) \varphi^2_n dx\right| \leq \|I_2\star w_\varepsilon^2\|_{L^\infty(\R^N)} \|\varphi_n\|^2_{L^2(Supp w_\varepsilon)} \rightarrow 0
\eeq
as $n \rightarrow \infty$. Here, we use the estimate: for any $x\in\R^N$
\beq\nonumber
\aligned
\left|I_2\star w_\varepsilon^2\right|(x)=&\int_{B_1(x)}\frac{w_\varepsilon(y)^2}{|x-y|^{N-2}}dy+\int_{B^c_1(x)}\frac{w_\varepsilon(y)^2}{|x-y|^{N-2}}dy \\
\leq &\|w_\varepsilon\|^2_{L^\infty} \int_{B_1(0)}\frac{1}{|y|^{N-2}}dy+\|w_\varepsilon\|_{L^2}^2<+\infty.
\endaligned
\eeq
We also have from the elementary inequality $|u_s|^2\leq 2|u_s-w_\varepsilon|^2+2|w_\varepsilon|^2$ that
$$\left|\int\left(I_2\star u^2_s\right) \varphi^2_n dx\right| \leq 2\|u_s-w_\varepsilon\|^{2}_{HL} \|\varphi_n\|^2_{HL} +2\int\left(I_2\star w_\varepsilon^2\right)  \varphi^2_n dx\leq \varepsilon^{2} \nu_s^{-1} +o(1).$$
This, \eqref{eq:20221121-9} and \eqref{eq:20221121-10} imply that
$$\int\left(I_2\star u_s^2\right)\varphi_n^2dx =o(1)~~~\text{as}~n \rightarrow \infty.$$
Moreover, by Corollary \ref{20221124-c1}
$$\int\left(I_2\star u_s\varphi_n \right) u_s \varphi_n dx=o(1).$$
Hence, by recalling Lemma \ref{3.1} we obtain
$$J''_s(u_s,\nu_s)[\varphi_n,\varphi_n]=\|\varphi_n\|^2_{H^s} -2\nu_s^2 \int\left(I_2\star u_s^2\right)\varphi_n^2dx-4\nu_s^2\int\left(I_2\star u_s\varphi_n \right) u_s \varphi_n dx=1+o(1).$$
This is in contradiction with \eqref{eq:20221121-7} and the proof of \eqref{eq:20221121-3} is complete.

(iii). Now we prove \eqref{eq:20221122-1}. Assume by contradiction that for every $s_0 \in (0,1)$ $$\inf\limits_{s \in \left(s_0,1\right]} \inf\limits_{u_s \in \mathcal{M}_s} K_r(s,u_s)=0.$$
Then there exist a sequence $s_n \nearrow 1$ and radial minimizers $u_{s_n}$ for $\nu_{s_n}$ such that
\beq\label{eq:20221122-2}
K_r(s_n,u_{s_n}) \rightarrow 0~~~\text{as}~n\rightarrow\infty.
\eeq
For fixed $n \in \N$, by the Riesz representation theorem and Ekeland variational principle, we obtain that there exist $f_{n,m} \in \Lambda^r_s$ and a minimizing sequence $\psi_{n,m} \in \Lambda^r_{s_n}$ for $K_r(s_n,u_{s_n})$ such that \beq\label{eq:20221225-4}
\|\psi_{n,m}\|_{H^{s_n}}=1~~~\forall m \in \N
\eeq
and
\beq\label{eq:20221122-4}
J''_s(u_{s_n},\nu_{s_n})[\psi_{n,m},v]-K_r(s_n,u_{s_n})\left\langle \psi_{n,m},v \right\rangle_{H^{s_n}}=\left\langle f_{n,m},v \right\rangle_{H^{s_n}},~~~\forall v \in \Lambda^r_{s_n},
\eeq
where $\|f_{n,m}\|_{H^{s_n}} \rightarrow 0$ as $m \rightarrow \infty$. Then there exists a sequence of sub-indices $m_n$ such that $\|f_{n,m_n}\|_{H^{s_n}} \rightarrow 0$ as $n \rightarrow \infty$.

Recalling Lemma \ref{3.5}, we may assume that $\nu_{s_n} \rightarrow \nu_1$ and, by Lemma \ref{20221103-l1}, that
\beq\label{eq:20221122-3}
\|u_{s_n}-U_1\|_{H^2} \rightarrow 0~~~\text{as}~n\rightarrow\infty.
\eeq
In particular, from \eqref{eq:20221122-4} we have
\beq\label{eq:20221122-5}
J''_s(u_{s_n},\nu_{s_n})[\psi_{n,m_n},v]-K_r(s_n,u_{s_n})\left\langle \psi_{n,m_n},v \right\rangle_{H^{s_n}}=\left\langle f_{n,m_n},v \right\rangle_{H^{s_n}}.
\eeq
\allowbreak
Let $w \in C^\infty_c (\R^N) \cap \Lambda^r_1.$ Then, from \eqref{eq:20221117-7} and \eqref{eq:20221122-3} we have
$$\left\langle \widehat{w},\widehat{u_{s_n}} \right\rangle_{H^{s_n}}=\int (1+|\xi|^{2s_n})\widehat{u_{s_n}} \widehat{w} d\xi= o(1)\|\widehat{w}\|_{H^{2+1}},$$
 and
$$\int (1+|\xi|^{2s_n})(i\xi^j)\widehat{u_{s_n}} \widehat{w} d\xi=o(1)\|\widehat{w}\|_{H^{2+2}},~~~\forall j=1,\cdots,N.$$
Let
$$v_n=w -\frac{\left\langle w,u_{s_n} \right\rangle_{H^{s_n}}}{\|u_{s_n}\|^2_{H^{s_n}}} u_{s_n}\in \Lambda^r_{s_n}.$$
 Using it as test function in \eqref{eq:20221122-5} and recalling that $\psi_{n,m_n} \in \Lambda^r_{s_n}$, we get
\beq\label{eq:20221122-6}
J''_s(u_{s_n},\nu_{s_n})[\psi_{n,m_n},w]-K_r(s_n,u_{s_n})\left\langle \psi_{n,m_n},w \right\rangle_{H^{s_n}}=o(1).
\eeq
\allowbreak
Since $\|\psi_{n,m_n}\|_{H^{s_n}}=1$, up to a subsequence, we have $\psi_{n,m_n} \rightharpoonup \psi~\text{in}~H^t(\R^N)$ for every fixed $t \in (0,1)$. Passing to the limit in \eqref{eq:20221122-6} and recalling \eqref{eq:20221122-2}, we get
$$J''(U_1,\nu_1)[\psi,w]=0~~~\forall w \in C^\infty_c(\R^N) \cap \Lambda^r_1.$$
Since, by Fatou's lemma, the latter identity implies that $\psi=0$, because the case $s=1$ is nondegenerate and $\psi \in \Lambda^r_1$.\par
That is, $\psi_{n,m_n} \rightharpoonup \psi=0$ in $H^t_r(\R^N)$ for every fixed $t \in (0,1)$ and so, by compactness, $$\psi_{n,m_n} \rightarrow 0~~~\text{in}~L^{\frac{4N}{N+2}}(\R^N).$$
Also, by \eqref{eq:20221122-5}, we have
\beq\label{eq:20221225-7}
J''_s(u_{s_n},\nu_{s_n})[\psi_{n,m_n},\psi_{n,m_n}]-K_r(s_n,u_{s_n}) \|\psi_{n,m_n}\|^2_{H^{s_n}}=o(1)
\eeq
and, by Corollary \ref{20221124-c1} and H-L-S inequality,
\beq\label{eq:20221122-7}
\aligned
\left|\int \left(I_2\star u_{s_n}\psi_{n,m_n}\right) u_{s_n}\psi_{n,m_n} dx\right|
\leq &\left|\int \left(I_2\star u^{2}_{s_n}\right) \psi^2_{n,m_n} dx\right| \\
 \leq &(\|u_{s_n}\|_{HL})^{2}(\|\psi_{n,m_n}\|_{HL})^{2}\\
 =&(\|\psi_{n,m_n}\|_{HL})^{2} \\
 \leq&\|\psi_{n,m_n}\|^2_{L^{\frac{4N}{N+2}}}=o(1),
\endaligned
\eeq
as $n \rightarrow +\infty$. Therefore, by \eqref{eq:20221225-7} and \eqref{eq:20221122-2} we get
\beq\nonumber
\aligned
&\|\psi_{n,m_n}\|^2_{H^{s_n}}-4\nu_{s_n}^2\int \left(I_2\star u_{s_n}\psi_{n,m_n}\right) u_{s_n}\psi_{n,m_n}dx-2\nu_{s_n}^2\int\left(I_2\star u^{2}_{s_n}\right) \psi^2_{n,m_n}dx  \\
=&o(1).
\endaligned
\eeq
Hence, passing to the limit, using \eqref{eq:20221225-4} and \eqref{eq:20221122-7}, we get $1-0=0$, that is a contradiction.
\ep

{\subsection{ Construction of pseudo-minimizers}

 Define the mapping
\beq\nonumber
\Phi_s : H^s_{rad}(\R^N ) \mapsto H^s_{rad}(\R^N ),\quad \Phi_s(\omega) = J'_s(U_1 + \omega, \nu_s ),
\eeq
$\Phi_s(\omega)$ and $J'_s(U_1 + \omega, \nu_s )$ are equal in the following sense:
for any $w \in  H^s_{rad}(\R^N )$, by Riesz representation theorem, there exists $\Phi_s(\omega)$ such that
\beq\nonumber
\langle\Phi_s(\omega), w\rangle_{H^s}=J'_s(U_1 + \omega, \nu_s )[w].
\eeq

\bl\label{4.2} For every $f \in H^s_{rad}(\R^N )$, there exists a unique $\bar{w}^s \in H^s_{rad}(\R^N )$ such that
\beq\nonumber
\langle\Phi'_s(0)[\bar{w}^s],w\rangle_{H^s}=\langle f,w\rangle_{H^s}\quad \forall w\in H^s_{rad}(\R^N ).
\eeq
In addition there exists a constant $C_1 > 0$ such that
\beq\label{eq:20221123-e1}
\|(\Phi'_s(0))^{-1}\|\leq C_1 \quad\forall s \in (s_0, 1).
\eeq
\el
\bp
We observe that $$\langle \Phi'_s(0)[w'],w\rangle_{H^s}=J''_s(U_1,\nu_s)[w',w].$$
Hence solving the equation $$\langle \Phi'_s(0)[\bar{w}],w\rangle_{H^s}= \langle f,w \rangle_{H^s}~~~\forall w \in H^s_{rad}(\R^N)$$
is equivalent to find a solution $\bar{w}$ to the equation
\beq\label{eq:20221123-1}
J''_s(U_1,\nu_s)[\bar{w},w]=\langle f,w \rangle_{H^s},
\eeq
for any $w \in H^s_{rad}(\R^N)$. By H-L-S inequality and H\"{o}lder inequality, for every $w \in H^s_{rad}(\R^N)$ we have
\begin{equation}\label{eq:20221123-2}
\aligned
&|(J''_s(U_1,\nu_s)-J''_s(u_s,\nu_s))[w,w]| \\
=&2\nu_s^2 \left|\int(I_2\star (u^{2}_s-U^{2}_1))w^2 dx+2\int(I_2\star (u_s - U_1) w)u_swdx \right. \\
&\left.+2\int(I_2\star U_1 w)(u_s-U_1)w dx\right| \\
\leq &2\nu_s^2 C\|u_s-U_1\|_{L^{\frac{4N}{N+2}}}\|w\|^2_{L^{\frac{4N}{N+2}}}.
\endaligned
\end{equation}
From Lemma \ref{20221103-l1} and Lemma \ref{3.5} we know that $\|u_s-U_1\|_{H^s} \rightarrow 0$ and $\nu_s \rightarrow \nu_1$ as $s \nearrow 1$. This implies that $u_s \rightarrow U_1$ in $L^{\frac{4N}{N+2}}(\R^N)$.
Therefore, from \eqref{eq:20221123-2},
\beq\label{eq:20221123-3}
|J''_s(U_1,\nu_s)-J''_s(u_s,\nu_s)[w,w]|=o(1)\|w\|^2_{H^s},
\eeq
where $o(1)$ is an infinitesimal quantity as $s \nearrow 1$.
This together with \eqref{eq:20221122-1} and \eqref{eq:20221121-2} in Lemma \ref{4.1} implies that there exist $C,s_0>0$ such that for all $s \in (s_0,1)$
\beq\label{eq:20221123-4}
|J''_s(U_1,\nu_s)[v,v]|=\left|J''_s(u_s,\nu_s)[v,v]+ o(1)\|v\|^2_{H^s}\right|\geq C\|v\|^2_{H^s}~~~\forall v \in H^s_{rad}(\R^N).
\eeq
\allowbreak
Hence, by the Lax-Milgram theorem, there exists a unique $\bar{w}^s \in H^s_{rad}(\R^N)$ such taht $$J''_s(U_1,\nu_s)[\bar{w}^s]=f$$
and by \eqref{eq:20221123-4}
$$\|\bar{w}^s\|_{H^s} \leq C\|f\|_{H^s},$$ which gives the desired result.
\ep

\bo\label{20221103-p1} For every $r > 0$ and $s > 0$, let
$$B_{r,s} = \left\{ \omega\in H^s_{rad}(\R^N ) : \|\omega\|_{H^s}\leq r \max\{1 - s, |\nu_1 - \nu_s |\}\right\}.$$
Then there exist $s_0 \in (0, 1)$ and $r_0 > 0$ such that for any $s \in (s_0, 1)$, there exists a unique
$\omega^s \in B_{r_0,s_0}$ such that
$\Phi_s(\omega^s ) = 0.$
\eo

\bp
We transform the equation $\Phi_s(\omega)=0$ to a fixed point problem:
\beq\label{eq:20221123-5}
\omega=-(\Phi'_s(0))^{-1}\{\Phi_s(0)+Q_s(\omega)\},
\eeq
\allowbreak
where $$Q_s(\omega):=\Phi_s(\omega)-\Phi_s(0)-\Phi'_s(0)[\omega].$$
$Q_s(\omega)$ is well-posed thanks to \eqref{eq:20221123-e1}. We observe that if $\omega \in H^s_{rad}(\R^N)$ then the mapping $$\omega \mapsto (\Phi'_s(0))^{-1}\{\Phi_s(0)+Q_s(\omega)\}$$ is radial too, since $U_1$ is radial.

For very $\bar{\omega} \in H^s_{rad}(\R^N)$, we set
\begin{align*}
\mathcal{N}_s(\omega)[\bar{\omega}]:=&\langle Q_s(\omega),\bar{\omega} \rangle_{H^s}\\
=&J'_s(U_1+\omega,\nu_s)[\bar{\omega}]-J'_s(U_1,\nu_s)[\bar{\omega}]-J''_s(U_1,\nu_s)[\omega,\bar{\omega}] \notag\\
=&2\nu_s^2\left(-\int(I_2\star (U_1+\omega)^2) (U_1+\omega)\bar{\omega}dx +\int(I_2\star (U_1)^2)U_1\bar{\omega}dx \right.  \\
&\left.+\int (I_2\star (U_1)^2) \omega \bar{\omega}dx+2\int (I_2\star U_1 \omega ) U_1 \bar{\omega}dx\right)\\
=&2\nu_s^2\left(-2\int(I_2\star (U_1\omega))\omega\bar{\omega}dx-\int(I_2\star \omega^2) (U_1+\omega)\bar{\omega}dx\right).\notag
\end{align*}
Then by H-L-S inequality, H\"{o}lder inequality and Sobolev inequality, we obtain
$$|\mathcal{N}_s(\omega)[\bar{\omega}]| \leq C(\|\omega\|^2_{H^s}+\|\omega\|^3_{H^s})\|\bar{\omega}\|_{H^s}$$ and
$$\|\mathcal{N}_s(\omega_1)-\mathcal{N}_s(\omega_2)\| \leq C(\|\omega_1\|_{H^s} +\|\omega_1\|^{2}_{H^s} +\|\omega_2\|_{H^s} +\|\omega_2\|^{2}_{H^s}) \|\omega_1-\omega_2\|_{H^s}.$$
This implies that for every $\|\omega_1\|_{H^s}$,$\|\omega_2\|_{H^s}<1$,
\beq\label{eq:20221123-7}
\|Q_s(\omega_1)\|_{H^s} =\|\mathcal{N}_s(\omega_1)\| \leq C_3\|\omega_1\|^{2}_{H^s}
\eeq
\allowbreak
and
\beq\label{eq:20221123-8}
\|Q_s(\omega_1)-Q_s(\omega_2)\|_{H^s}=\|\mathcal{N}_s(\omega_1)-\mathcal{N}_s(\omega_2)\| \leq C_3\|\omega_1-\omega_2\|_{H^s},
\eeq
\allowbreak
where $C_3$ is independent on $s \in (s_0,1)$.

Now we claim that there exists a constant $C_2>0$ independent on $s \in (s_0,1)$ such that
\beq\label{eq:20221123-9}
\|\Phi_s(0)\| \leq C_2 \max \{1-s,|\nu_1 -\nu_s|\}.
\eeq
\allowbreak
Indeed, we have from \eqref{eq:20221117-7} that
$$|J'_s(U_1,\nu_s)[v]-J'_1(U_1,\nu_1)[v]| \leq (1-s)C_{\delta,N}\|U_1\|_{H^{2-s+2\delta}}\|v\|_{H^s} +|\nu_1-\nu_s|\|v\|_{H^s}.$$
Then by \eqref{eq:20220913-e3}, $J'_1(U_1,\nu_1)=0$, we get \eqref{eq:20221123-9}.

Now we solve the fixed point equation \eqref{eq:20221123-5} in the ball
$$B_{r,s}=\{\omega \in H^s_{rad}(\R^N): \|\omega\|_{H^s} \leq r\alpha_s\}$$
to finish the proof of Proposition \ref{20221103-p1},
where $\alpha_s=\max\{1-s,|\nu_1-\nu_s|\}$ and $r>0$ will be selected later. Indeed for $\omega \in B_{r,s}$, by \eqref{eq:20221123-e1}, \eqref{eq:20221123-9} and \eqref{eq:20221123-7}, we have
$$\|\Phi'_s(0)^{-1}\{\Phi_s(0) +Q_s(\omega)\}\|_{H^s}\leq C_1(C_2 \alpha_s+C_3 r^2 \alpha^2_s).$$
There exists $r_0>0$ large and $s_0=s_0(r_0) \in (0,1)$ such that for any $s \in (s_0,1)$ we have
$$C_1(C_2\alpha_s +C_3r^2 \alpha^2_s)<C_1(C_2\alpha_{s_0} +C_3r^2_0 \alpha^2_{s_0})\leq r_0\alpha_{s_0}.$$
Since $\lim\limits_{s\nearrow 1}\alpha_s=0$, then for every $s \in (s_0,1)$ the mapping
$$\omega \mapsto -(\Phi'_s(0))^{-1}\{\Phi_s(0)+Q_s(\omega)\}$$
maps $B_{r_0,s_0}$ into itself. Select a suitable $s_0$, this map is a contraction on $B_{r_0,s_0}$ by \eqref{eq:20221123-8}.

Therefore, by the Banach fixed point theorem, for every $s \in (s_0,1)$, there exists a unique function $\omega^s \in B_{r_0,s_0}$ solving the fixed point equation \eqref{eq:20221123-5}.

\ep

{\subsection{  Completion of the proof of Theorem \ref{th1.3}. }

By using the same argument of Sec. 5.2 in \cite{FV}, we get the uniqueness result.

\vskip4mm
%%%%%%%%%%%%%%%%%%%%

\end{document}